# Non-Markovian superposition process model for stochastically describing concentration–discharge relationship

Hidekazu Yoshioka[a,*] and Yumi Yoshioka[b]

[a]Japan Advanced Institute of Science and Technology, 1-1 Asahidai, Nomi 923-1292, Japan
[b]Gifu University, 1-1 Yanagido, Gifu 501-1193, Japan
[*]Corresponding author: yoshih@jaist.ac.jp, ORCID: 0000-0002-5293-3246

**Abstract**
Concentration–discharge relationship is crucial in river hydrology, as it reflects water quality dynamics across both low- and high-flow regimes. However, its mathematical description is still challenging owing to the underlying complex physics and chemistry. This study proposes an infinite-dimensional stochastic differential equation model that effectively describes the concentration–discharge relationship while staying analytically tractable, along with the computational aspects of the model. The proposed model is based on the superposition of the square-root processes (or Cox–Ingersoll–Ross processes) and its variants, through which both the long-term moments and autocovariance of river discharge and the fluctuation of water quality index can be derived in closed forms. Particularly, the model captures both long (power decay) and short (exponential decay) memories of the fluctuation in a unified manner, while quantifying the hysteresis in the concentration–discharge relationship through mutual covariances with time lags. Based on a verified numerical method, the model is computationally applied to weekly data on total nitrogen (TN. long memory with moderate fluctuation), total phosphorus (TP. short memory with large fluctuation), and total organic carbon (TOC. short memory with moderate fluctuation) from a rural catchment to validate its applicability to real-world datasets. Based on the identified model and its mutual covariance, our findings indicate that, on average, the peak concentrations of these water quality indices appear approximately 1 day after discharge. Finally, the study discusses the effects of model uncertainty on mutual covariance.

**Statements & Declarations**
**Acknowledgments:** The authors would like to express their gratitude to Dr. Ikuo Takeda of Shimane University for providing water quality data at Osu.
**Funding:** This study was supported by the Japan Society for the Promotion of Science (KAKENHI No. 22K14441) and Japan Science and Technology Agency (PRESTO No. JPMJPR24KE).
**Competing Interests:** The authors have no relevant financial or non-financial interests to disclose.
**Data Availability:** The data will be made available upon reasonable request to the corresponding author.
**Declaration of Generative AI in Scientific Writing:** The authors did not use generative AI for the scientific writing of this manuscript.

## 1. Introduction

### 1.1 Research Background

Rivers function as conduits for water, organic and inorganic materials, and aquatic species, while also providing water resources for diverse human activities, including irrigation, drinking water supply, and hydropower generation (Andrade-Muñoz et al. 2025; Cowx et al. 2024; Jiang et al. 2024; Juárez et al. 2023)[1,2,3,4]. Consequently, studying the quantity and quality of river water has become a pivotal focus in hydrology and related research fields given its critical role in environmental risk and uncertainty assessments (Chaudhary et al. 2024; Elagib et al. 2021; Srinivasan et al. 2021)[5,6,7].

River-water quantity (streamflow discharge) and quality (water quality) are interdependent, with the former influencing the latter through short- and long-term dynamics, collectively termed concentration–discharge relationships (Jing et al. 2025; Pohle et al. 2021)[8,9]. Over longer timescales (daily, monthly, or beyond), correlations between streamflow discharge and water quality indices reveal catchment characteristics such as land use, geology, and climate (da Conceição et al. 2024; Diaz et al. 2024; Kaltenecker et al. 2025; Stenger et al. 2024)[10,11,12,13]. Conversely, at shorter timescales (hourly to daily) and event scales, concentration–discharge curves capture the inertia or hysteresis in water quality dynamics relative to the discharge fluctuations during flood events (Liu et al. 2021; Mazilamani et al. 2024; Musolff et al. 2024)[14,15,16]. The concentration–discharge relationship aids in understanding hydrological processes within a catchment, such as hydrograph separation (Woodward and Stenger 2018; Woodward and Stenger 2020)[17,18]. Von Brömssen et al. (2023)[19] employed a nonparametric method to analyze the concentration–discharge relationship that would not be represented by a simple function.

Event-based concentration–discharge relationships are often quantified using hysteresis indices, which, given time-series data, involve the following: (1) Identifying a flood event, (2) plotting the corresponding concentration–discharge loop in a two-dimensional phase space defined by concentration and discharge, (3) distinguishing the rising and falling limbs of the loop, and (4) calculating specific integrals along the curve (Zuecco et al., 2016; Gao et al. 2024; Heathwaite and Bieroza 2016; Li et al. 2022; Lloyd et al. 2016; Mehdi et al. 2021)[20,21,22,23,24,25]. These methods are effective when a loop and its limbs can be easily identified; however, identification is often challenging because concentration–discharge loops during a single flood event exhibit multiple crossings, such as figure-eight shapes, or even more complicated patterns (Ariano and Ali 2024; Lannergård et al. 2021; Yue et al. 2023)[26,27,28]. Therefore, a unified method for evaluating hysteresis in the concentration–discharge relationship is necessary.

The stochastic nature of surface-water quantity and quality, which fluctuate randomly over time and often exhibit long memories characterized by subexponentially decaying autocorrelation with fractional dynamics (Liang et al. 2021; Sun et al. 2024; Yu et al. 2023; Pärn et al. 2024; Zhang et al., 2018)[29,30,31,32,33], poses additional challenges in analyzing concentration–discharge relationships. This behavior suggests the use of some non-Markovian, memory-dependent processes (Yoshioka and Yoshioka 2024a)[34]. A mathematical model that provides a unified description of multi-scale water quantity and quality dynamics with a suitable memory structure is valuable for both scientific and engineering purposes. Such a unified model can be applied to various problems, including, but not limited to, the statistical analysis of hydrological processes and impact studies of human activities on catchment hydrology. Accordingly, model interpretability is important because increasing model complexity, particularly in a nonparametric pattern such as statistical and machine learning approaches, can improve model performance but may largely depend on unexplainable causality relationships. Moreover, analytically tractable models would become powerful tool in applications if their parameters and coefficients are identified using common statistical methods such as least squares and moment matching.

Theoretical and technical challenges must be addressed to establish a versatile stochastic process model capable of accommodating multiple timescales, analytical tractability, and interpretability. Key challenges include the non-Markovian nature resulting from long memory and computational complexities, some of which have been partly resolved. For memory, a non-Markovian process can be formally decomposed into a sum or a proper integral of Markovian processes, referred to as Markovian lifts or Markovian embedding (Cuchiero and Teichmann 2020; Kanazawa and Sornette 2024)[35,36]. The approach enables the handling of simpler, higher-dimensional processes that correspond to difficult lower-dimensional dynamics. Within hydrology, Markovian lifts decompose processes with long memories into numerous smaller processes running on distinct timescales. This approach has been effective in modeling streamflow discharge (Yoshioka and Yoshioka 2024b)[37] and water quality indices (Yoshioka and Yoshioka 2024a)[34], although their interactions have not been considered. Yoshioka and Yoshioka (2025)[38] recently demonstrated that applying a variant of the stochastic volatility model, commonly used in finance and economics (Todorov, 2011)[39], could effectively represent the "logarithm" of water quantity dynamics. However, mathematical instability issue caused by the specific structure of the model limits its usefulness.

Computability presents a significant challenge, particularly in the modeling and computation of water quality dynamics. Typically, the water quality index represents the concentration of chemicals in the material and must remain non-negative. Non-negativity can be theoretically preserved by modeling water quality dynamics using specific stochastic differential equations (SDEs), particularly the square-root process (called Cox–Ingersoll–Ross process) and its variants (Song et al. 2023; Wang et al. 2024; Zonouz et al. 2024)[40,41,42]. These processes are appealing for applications because their statistics, including cumulants and autocorrelation, can be derived in closed forms (e.g., Alfonsi, 2015)[43]. However, generating their sample paths—simulating them—proved to be extremely challenging, as highlighted in the literature on theoretical numerical analysis. Common numerical methods for SDEs, such as the Euler–Maruyama method, exhibit arbitrarily slow convergence rates when diffusion dominates (Hefter and Herzwurm 2018; Mickel and Neuenkirch 2025)[44,45]. Thus, reliable simulation is often unattainable unless specialized discretization schemes, which are not always easy to implement, are used (Multi-level Monte–Carlo method (Zheng 2023)[46] and machine learning techniques (van Rhijn et al. 2023)[47]). Accurate numerical methods exist but are limited to cases with small diffusion (Kamrani and Hausenblas 2025; Yi et al. 2021)[48,49] or computing conditional expectations rather than generating sample paths (Alfonsi and Lombardo 2024; Alfonsi, 2025)[50,51]. This limitation arises from the non-linear and non-Lipschitz diffusion coefficients, which preserve the non-negativity of the solutions. Specifically, with a small diffusion, the square-root process provides positive values, while with a larger diffusion, it provides only non-negative values that may touch zero. This behavior was recently analyzed in detail by Mishura et al. (2024)[52] using path-wise analysis. The computability issue has hindered progress in understanding the non-Markovian processes governing water quantity and quality dynamics, as simulations are essential for analyzing the concentration–discharge curve during flood events. Additionally, the time multi-scale nature of water-quality dynamics inherently involves square-root diffusion with large diffusion coefficients (**Sections 2** and **3**) that are difficult to simulate using common numerical methods.

### 1.2 Objectives and Contributions

The objective of this study is to propose and examine a tractable, interpretable, and computable non-Markovian stochastic process model that can capture memory structures in the quantity and quality dynamics of river water and the concentration–discharge relationship. Our contributions to the research objective are as follows:

- The proposed model comprises a long-memory streamflow discharge process, generalizing the superposition of Ornstein–Uhlenbeck processes as the simplest non-Markovian process model (Barndorff-Nielsen 2001)[53]. We assume that the autocorrelation of our model is given by integrating exponential functions against a probability measure, such that Markovian lifts apply while keeping analytical tractability during statistical analysis (Fasen 2009; Fuchs et al. 2013)[54,55]. Consequently, the discharge of the proposed model builds on existing frameworks with certain generalizations, while the water quality aspect relies on the superposition of square-root processes, dependent on streamflow discharge–a novel modeling strategy. Dupret et al. (2023)[56] explored the superposition of square-root processes for finance in the stochastic Volterra framework.
- The resulting concentration–discharge model is a non-linear yet affine SDE (Duffie et al., 2003)[57], where cumulants are derived analytically or by solving a generalized Riccati equation, despite its non-Markovian nature (Bondi et al. 2024; Motte and Hainaut 2024)[58,59]; we employ the former approach. Each term in the model can be interpreted owing to its simple form. Additionally, autocorrelation and mutual covariance can be derived in closed forms; the former captures memory structures in the concentration–discharge relationship, while the latter assesses event-wise hysteresis between water quantity and quality from a statistical perspective without introducing external methodologies and indices.
- We demonstrate that the analytical tractability facilitates computational efficiency, despite being based on square-root processes, marking another contribution of this study. Further, we show that the dynamics-based discretization method proposed by Abi Jaber (2024)[60] can be applied to our model with a proper enhancement. This discretization method has been applied to square-root processes in finance and economics, but it has not been used in water quantity and quality dynamics.

A potential computational challenge in our context is that the stochastic components become significantly larger than the deterministic components in Markovian lifts, which precludes the use of common numerical discretization methods, as in the classical square-root process. The method of Abi Jaber (2024)[60] is still applicable to such cases, capturing cumulants and first hitting times—both statistical and pathwise quantities—thus supporting its application to fields beyond finance and economics, and to the concentration–discharge relationship. With this contribution, it would become possible to approach a broader class of stochastic models computationally.

- The final contribution of this study is the application of our model to time-series data of total nitrogen (TN), total phosphorus (TP), and total organic carbon (TOC), which were sampled weekly at an observation point in a river in Japan. The model parameters are identified by combining least-squares and moment-matching methods, leveraging their analytical tractability. Our findings indicate that, at the study site, discharge and TN include long memory, while TP and TOC display short memories. These water quality indices are positively correlated with discharge, and mutual covariances suggest that their highest concentrations occur approximately one day after the discharge peak. The proposed numerical method reasonably reproduces statistics and probability densities of the water quality indices. Computed sample paths of TN and TOC differ qualitatively from those of TP due to the relatively high noise intensity in TP; the TP paths resemble square-root processes with large noise intensity, often approaching zero. These refined statistical descriptions of the concentration–discharge relationship, based on an interpretable stochastic process model, are absent in the literature. In addition, any model fitted with data is subject to errors due to data quantity and quality, and the limitation of fitting methods, unless the model itself is the ground truth (Aibaidula et al., 2023; Bacci et al., 2023; Baig et al., 2025)[61,62,63]. We thus evaluate the influence of model uncertainty on the peak location in mutual covariance when there is a model misspecification.

Overall, this study provides a novel mathematical model for the concentration–discharge relationship, along with numerical computations and applications.

### 1.3 Structure of This Paper

The remainder of this paper is organized as follows. **Section 2** presents the proposed mathematical model and its statistics. **Section 3** focuses on a numerical method for simulating sample paths in the proposed model. **Section 4** describes model applications to field data. **Section 5** presents the conclusions drawn from the study. **Appendices** provide proofs (**Section A1**) and auxiliary results (**Sections A2** and **A3**).

## 2. Mathematical Model

The proposed model comprises a long-memory process of streamflow discharge at an observation point in a river and the superposition of square-root process and its variants driven by discharge that governs the (de-seasonalized) water quality index. The coupling direction is from discharge to concentration, based on the physical consideration that water quality dynamics in a river environment are shaped by runoff from the catchment to the river and that the movement of water is not affected by its quality. We consider a complete probability space characterized by triplet function $(\Omega, \mathbb{F}, \mathbb{P})$, denoting a collection of events $\Omega$, filtration $\mathbb{F}$ (generated by the white noise $\zeta$ and discharge $Q$ (or its driving noise process) explained later), and probability $\mathbb{P}$, as is typical in stochastic models (Chapter 1 by Alfonsi (2015)[43]). Expectation is represented as $\mathbb{E}$ and variance as $\mathbb{V}$. Time is considered a real parameter.

## 2.1 Discharge

Following the literature on hydrology (Botter et al. 2013; Müller et al. 2021; Miniussi et al. 2023)[64,65,66], we assume that the streamflow discharge $Q = (Q_t)_{t \in \mathbb{R}}$, which is defined as the water volume that passes through a river cross-section in unit time, is a continuous-time scalar and non-negative process driven by jumps. This assumption can be justified for long-term (e.g., longer time scales than single flood events) data having hourly to daily resolution considered in our application. In the model, a jump event conceptually represents a flood event. Non-negativity means no flow reversal.

We assume the discharge $Q$ to be stationary, with stationary average $\bar{Q} = \mathbb{E}[Q_t] > 0$ and variance $\bar{V} = \mathbb{V}[Q_t] > 0$. Moreover, we assume that the autocorrelation $A_Q(h)$ of $Q$ with lag $h \geq 0$ has the integral representation:

$$A_Q(h) = \int_0^{+\infty} e^{-rh} \theta(\mathrm{d}r) \tag{1}$$

with a probability measure $\theta$ for a positive random variable. Eq. (1) implies that the autocorrelation of discharge is a superposition (integration) of the exponential $e^{-rh}$ with respect to $\theta(\mathrm{d}r)$. Choosing a Dirac delta as $\theta$ reduces Eq. (1) to an exponential function corresponding to short memory. A long memory function is defined when $\theta$ has a density $p$, which is singular near $r = 0$ with the scaling relationship $p(r) = O(r^{c-1})$ with small $c > 0$ (Barndorff-Nielsen and Leonenko 2005)[67]. The integral representation (1) is limited by a decreasing and convex autocorrelation $A_Q$ owing to its convex combination of decreasing exponential functions. This drawback can be overcome by increasing the degrees of freedom of the in the integrand such as by adding sinusoidal functions (Yoshioka 2024)[68] but is not employed in this study for calculational simplicity. Such an approach is necessary when considering discharge in an estuary where the ebb and flow of the tide are non-negligible. A specific model for discharge is discussed in **Section 4**.

## 2.2 Concentration

We assume that the concentration of a water quality index $C = (C_t)_{t \in \mathbb{R}}$ is non-negative, and is driven by discharge $Q$ and fluctuations independent of it. We also assume that $C$ can be decomposed into non-negative seasonal $S$ and stationary stochastic components $M$, the latter being non-dimensional:

$$C_t = S_t M_t. \tag{2}$$

The seasonal component $S = (S_t)_{t \in \mathbb{R}}$ is deterministic with a period of one year. The stochastic component $M = (M_t)_{t \in \mathbb{R}}$, hereafter referred to as the memory process because the memory of water quality is encoded in it, is modeled by the superposition approach. This is possibly the simplest model that accounts for seasonality.

We formulate $M$ to satisfy three requirements.

- $M$ is non-negative and has correlation to discharge $Q$ as physically required.
- $M$ is an affine process so that the superposition becomes possible, and moments are obtained in closed forms.
- $M$ can address both short and long memories as observed in real data.

For the first and second requirements, we employ the square-root process as an affine process that is convenient to analyze due to explicit availability of moments. Moreover, an affine process is able to describe a non-negative stochastic process. So, it would be a candidate for modelling non-negative variables like water quality indices. The third requirement excludes the direct use of the classical square-root process that has an exponential autocorrelation. The superposition approach (Yoshioka and Yoshioka 2024a)[34] is extended to the proposed mathematical modeling, such that $M$ is discharge-dependent, affine, and can reproduce a wide range of memories.

Considering the required properties, $M$ is formulated as follows:

$$M_t = \underbrace{\int_0^{+\infty} m_t(\mathrm{d}R)}_{\text{Superposition of processes of square-root type}}, \quad t \in \mathbb{R}, \tag{3}$$

where $m = (m_t(\mathrm{d}R))_{t \in \mathbb{R}}$ is a measure-valued process that satisfies the square-root type SDE for each $R > 0$:

$$\mathrm{d}m_t(\mathrm{d}R) = -\underbrace{R}_{\text{Mean reversion}} \left( m_t(\mathrm{d}R) - \underbrace{(a + bQ_t)}_{\text{Discharge dependence}} \rho(\mathrm{d}R) \right) \mathrm{d}t + \underbrace{\sigma \sqrt{R} B(\mathrm{d}R, \mathrm{d}t)}_{\text{Diffusion}}, \tag{4}$$

where $a, b \geq 0$ are parameters to represent the discharge dependence of concentration and $\sigma \geq 0$ is noise intensity. Here, $B$ denotes space-time (and double-sided in time) Gaussian random measure with the formal state-dependent covariance (0, if $R \neq R'$ or $t \neq t'$):

$$\begin{aligned}
\mathbb{E}[B(\mathrm{d}R, \mathrm{d}t) B(\mathrm{d}R', \mathrm{d}t')] &= \mathbb{E}\left[\sqrt{m_t(\mathrm{d}R) m_{t'}(\mathrm{d}R')}\right] \mathbb{E}[\zeta(\mathrm{d}R, \mathrm{d}t) \zeta(\mathrm{d}R', \mathrm{d}t')] \\
&= \delta(R - R') \delta(t - t') \mathbb{E}\left[\sqrt{m_t(\mathrm{d}R) m_{t'}(\mathrm{d}R')}\right] \mathrm{d}t \mathrm{d}t'
\end{aligned}. \tag{5}$$

Here, $\delta$ is the Dirac delta, $\rho(\mathrm{d}R)$ is the probability measure of reversion $R > 0$, $\zeta$ is the white noise

(Chapter 3 by Pardoux (2021))[69]), $a, b \geq 0$ are parameters to represent the discharge dependence of concentration, and $\sigma \geq 0$ is noise intensity. We assume that $\zeta$ is independent of discharge $Q$. Eq. (4) is formally solved as follows by an integration:

$$m_t(\mathrm{d}R) = \int_{s=-\infty}^{s=t} e^{-R(t-s)} \left\{ R(a+bQ_s)\rho(\mathrm{d}R)\mathrm{d}s + \sigma\sqrt{R}B(\mathrm{d}R,\mathrm{d}s) \right\}. \quad (6)$$

In Eq. (3), the integration with respect to $R$ suggests that the concentration time series contains multiple time scales, as $R$ can physically represent the inverse of timescale for each $m$. The mean reversion term in (4) means that the process $m_t(\mathrm{d}R)$ reverts to its average with the time scale $R^{-1}$. The inclusion of the probability measure $\rho$ in the drift of (4) is because $m_t(\mathrm{d}R)$ is measure-valued.

The proposed model is non-Markovian because the conditional expectation of the concentration $C_t$ at time $t$ cannot be predictable based on some past value $C_s$ with $s < t$, but with the full information at time $s$ (i.e., $\mathbb{F}_s$ or whole $m_s(\mathrm{d}R)$); collecting the latter information needs observing the whole history of concentration $C$.

The process $m$ is formally expressed as the square-root process

$$\mathrm{d}m_t = -R\left(m_t - (a+bQ_t)\rho\right)\mathrm{d}t + \sigma\sqrt{Rm_t}\,\mathrm{d}W_t, \quad (7)$$

where $W$ is a one-dimensional standard Brownian motion. The representation (7) is ambiguous because $\rho$ is not a function but a probability measure. Nevertheless, Eq. (7) is useful for understanding the origin of the integration with respect to $R$ performed in Eq. (6) from a finite-dimensional perspective where the probability measure $\rho$ can be replaced by an empirical measure (**Section A3**). This understanding is important for analyzing the computational aspects of the proposed model (**Section 3**).

Finally, from Eqs. (2), (3) and (6), the concentration $C$ is obtained as follows:

$$C_t = S_t \int_{R=0}^{R=+\infty} \int_{s=-\infty}^{s=t} e^{-R(t-s)} \left\{ R(a+bQ_s)\rho(\mathrm{d}R)\mathrm{d}s + \sigma\sqrt{R}B(\mathrm{d}R,\mathrm{d}s) \right\}. \quad (8)$$

We focus on the analysis of memory process $M$ in the next subsection because $S$ is deterministic.

***Remark 1*** One may consider incorporating seasonality into the coefficients of SDEs (Mostafa and Allen 2024)[70]. In our framework, the consideration would be possible at the cost of a significant loss of analytical tractability.

***Remark 2*** One may consider simpler non-negative processes such as geometric Brownian motion (Taylor et al. 2022)[71] for superposition. However, these models are not affine, and the superposition approach has not yet been examined. Conversely, the superposition of affine processes, including square-root processes and their variants, has been systematically studied (Barndorff-Nielsen 2001; Iglói and Terdik 2003)[53,72]. Therefore, we use the SDE (4). The analytical tractability of our model is attributable to the square-root coefficient ($\sqrt{Rm_t}$ in Eq. (7)), using which we formally derive the averages of $\mathbb{E}[m_t] = O(\mathrm{d}R)$ and $\mathbb{V}[m_t] = O(\mathrm{d}R)$ as $\mathbb{E}[M] = O(1)$ and $\mathbb{V}[M_t] = O(1)$, respectively. These relationships suggest that $M$ has finite average and variance.

## 2.3 Statistics

### 2.3.1 Average and variance

The average and variance of the memory process $M$ are obtained analytically by considering the expectations of Eq. (3). The proof is based on a formal calculation that can also be deduced from a finite-dimensional version of the system (3)-(4); see **Section A2 in Appendix** for the finite-dimensional system.

***Proposition 1***

*Assuming a stationary state, it follows that*

$$\mathbb{E}[M_t] = a + b\bar{Q} \tag{9}$$

*and*

$$\mathbb{V}[M_t] = \underbrace{\frac{1}{2}\sigma^2 \mathbb{E}[M_t]}_{\text{Fluctuation not directly coming from discharge}} + \underbrace{b^2 \bar{V} \int_{R=0}^{R=+\infty} \int_{P=0}^{P=+\infty} \int_{r=0}^{r=+\infty} \frac{PR}{P+R}\left(\frac{1}{P+r} + \frac{1}{R+r}\right) \theta(\mathrm{d}r)\rho(\mathrm{d}R)\rho(\mathrm{d}P)}_{\text{Fluctuation directly coming from discharge}}. \tag{10}$$

According to **Proposition 1**, the average of the memory process $M$ depends linearly on the average of discharge, which can be considered as the simplest dependence between water quality and discharge. Specifically, the proposed model assumes concentration to increase with discharge in the long run, implying that the model covers cases where a positive correlation exists between concentration and discharge. **Proposition 3** demonstrates that the proposed model has a positive correlation. Higher order cumulants such as skewness and kurtosis may be available, but with integrals having higher dimensions.

### 2.3.2 Autocorrelation

The autocorrelation $A_M(h)$ of $M$ with lag $h \geq 0$ is obtained in a closed form.

***Proposition 2***

*Assuming a stationary state, the autocorrelation $A_M(h)$ ($h \geq 0$) is derived in a closed form as follows:*

$$A_M(h) = \frac{1}{\mathbb{V}[M_t]} \left( \underbrace{\frac{\sigma^2}{2} \mathbb{E}[M_t] \int_0^{+\infty} e^{-Rh} \rho(\mathrm{d}R)}_{\text{Correlation not directly coming from discharge}} + \underbrace{b^2 \bar{V}\left(I(h) + J(h)\right)}_{\text{Correlation directly coming from discharge}} \right) \tag{11}$$

*with*

$$I(h) = \int_{R=0}^{R=+\infty} \int_{P=0}^{P=+\infty} \int_{r=0}^{r=+\infty} \frac{RP}{P+R}\left(\frac{1}{P+r} + \frac{1}{R+r}\right) e^{-Ph} \theta(\mathrm{d}r)\rho(\mathrm{d}R)\rho(\mathrm{d}P) \tag{12}$$

*and*

$$J(h) = \int_{R=0}^{R=+\infty} \int_{P=0}^{P=+\infty} \int_{r=0}^{r=+\infty} \frac{RP}{(R+r)(P-r)} \left(e^{-rh} - e^{-Ph}\right) \theta(\mathrm{d}r) \rho(\mathrm{d}R) \rho(\mathrm{d}P). \tag{13}$$

According to **Proposition 2**, the autocorrelation $A_M$ consists of two parts: fluctuation unrelated to discharge and that related to discharge. The asymptotic behavior of each integral for a large $h > 0$ depends on the regularity of the probability measures $\theta$ and $\rho$ near the origin. The first integral in Eq. (11) and $I(h)$ strictly decrease with respect to $h > 0$, whereas $J(h)$ may be unimodal with a strict maximum point at some $h > 0$.

To investigate the existence of autocorrelation and variance of $M$, we apply Proposition 3 in Yoshioka and Yoshioka (2025)[38] to our **Propositions 1–2** owing to similarities in functional shapes, despite the models being fundamentally different. Specifically, following this literature, we assume that $\rho$ has a gamma distribution that behaves asymptotically near $R = 0$ as $O(R^{\alpha_R - 1})$ ($\alpha_R > 0$). Further, we consider $\theta$ as another gamma density that behaves asymptotically near $r = 0$ as $O(r^{\alpha_r - 1})$ ($\alpha_r > 0$), as in **Section 4**. Here, variance and autocorrelation exist if $\alpha_r > 1$, or if $\alpha_r \in (0,1]$ and $\alpha_r(\alpha_R + 2) - 1 > 0$, such that the underlying memories of the driving noise processes are not too large. This sufficient condition holds true in our case study in **Section 4**.

### 2.3.3 Mutual Covariance

In this subsection, the dynamic correlation between $C$ and $Q$ is analyzed through the mutual covariance between $M$ and $Q$, as the latter can derived in a closed form. The theoretical assumption of annually varying seasonal factor $S$ along with the physical assumption that each flood event lasts on a daily to weekly scale, as satisfied in our application in **Section 4**, suggests that the mutual covariance between $C$ and $Q$ can be effectively analyzed through that of $M$ and $Q$.

With time lag $h \in \mathbb{R}$, the mutual covariance is set as follows:

$$\lambda(h) = \mathbb{E}\left[(Q_t - \mathbb{E}[Q_t])(M_{t+h} - \mathbb{E}[M_t])\right]. \tag{14}$$

A positive $h$ (resp., negative $h$ case) evaluates the correlation between future memory process $M_{t+h}$ and past discharge $Q_t$ (resp., past memory process $M_{t+h}$ and future discharge $Q_t$). Specifically, the mutual covariance $\lambda$ evaluates the delay of peaks between $M$ and $Q$. The event-scale concentration–discharge relationship can be characterized by a loop in the two-dimensional $C$-$Q$ phase space (Mazilamani et al. (2024))[15], where the shape and direction of the loop can be determined by the sign and size of lag (the sign and absolute value of $h$) between concentration and discharge peaks. Therefore, $\lambda$ effectively characterizes such loops, and $\lambda$ can be obtained in a closed form with suitable integrals as shown below.

***Proposition 3***

*Assuming that $\int_0^{+\infty} R\rho(\mathrm{d}R) < +\infty$, the mutual covariance $\lambda(h)$ ($h \in \mathbb{R}$) at a stationary state is derived in a closed form as follows:*

$$\lambda(h) = \begin{cases} F(h) & (h \geq 0) \\ G(-h) & (h < 0) \end{cases} \quad (15)$$

with

$$F(l) = b\bar{V}\int_{R=0}^{R=+\infty}\int_{r=0}^{r=+\infty}\left\{\frac{R}{R+r}e^{-Rl} + \frac{R}{R-r}\left(e^{-rl} - e^{-Rl}\right)\right\}\theta(\mathrm{d}r)\rho(\mathrm{d}R), \ l \geq 0 \quad (16)$$

and

$$G(l) = b\bar{V}\int_{R=0}^{R=+\infty}\int_{r=0}^{r=+\infty}\frac{R}{R+r}e^{-rl}\theta(\mathrm{d}r)\rho(\mathrm{d}R), \ l \geq 0. \quad (17)$$

*As a byproduct, the covariance between $M$ and $Q$ is derived as follows:*

$$\mathbb{E}\left[(Q_t - \mathbb{E}[Q_t])(M_t - \mathbb{E}[M_t])\right] = \lambda(0) = b\bar{V}\int_{R=0}^{R=+\infty}\int_{r=0}^{r=+\infty}\frac{R}{R+r}\theta(\mathrm{d}r)\rho(\mathrm{d}R). \quad (18)$$

The assumption $\int_0^{+\infty} R\rho(\mathrm{d}R) < +\infty$ in **Proposition 3** means that the probability measure $\rho$ has an average, which is not restrictive in applications.

According to **Proposition 3**, $\lambda$ is continuous at $h = 0$, with positive and negative components in Eqs. (16) and (17). The negative component (Eq. (17)), which determines the delay of $Q$ relative to $M$, is strictly increasing in $h$, suggesting that the correlation between future discharge and past concentration increases as the delay decreases. Conversely, the positive component (Eq. (17)), which evaluates the delay of $M$ relative to $Q$, is not necessarily monotonic with respect to $h$. As discussed in **Section 4**, the models identified based on real datasets suggest that Eq. (16) exhibits one maximum point at some $h = h_c > 0$. Indeed, the integrand in Eq. (16) is the sum of a decreasing function and a unimodal function of $h$, and the summed functions possibly exhibit a strict maximum at some $h > 0$ when the latter function is sufficiently large.

The maximum of the positive component at $h = h_c > 0$, when it exists, corresponds to the peak delay between the concentration (more rigorously, the memory process) and the past discharge. Herein, the form of $\lambda$ suggests that, on average, the peak concentration appears after $h_c$ from that of the discharge, implying the dominance of counter-clockwise hysteresis of the concentration–discharge loop often reported in the literature (Liu et al. 2021; Hou et al. 2023; Marin-Ramirez et al. 2024)[73,74,75]. This type of hysteresis is consistent with the physical assumption considered for the proposed model, particularly the drift of the SDE (4), where concentration is partly driven by discharge in a positive manner, such that concentration increases with discharge. However, the assumption does not imply that the proposed model generates only the counterclockwise hysteresis of the concentration–discharge loop. Clockwise-type hysteresis concentration–discharge loops can be observed when fluctuations are not directly related to discharge in the SDE (4), i.e., the last term in Eq. (4).

The proposed model poses limitations such that it cannot capture the negative correlation between concentration and discharge when observed for natrium ions (Floury et al. 2024)[76], considering specific conductance (Clark et al. 2019)[77]. Although this limitation can be overcome by simply replacing the coefficient $a+bQ_t$ by $a-bQ_t$ or $-a-bQ_t$; however, replacement cannot be performed in the framework of affine processes because the process $m$ possibly becomes negative as $Q_t$ is theoretically unbounded from above. Accordingly, the stochastic volatility model (Yoshioka and Yoshioka, 2025)[38] is more advantageous although it considers the "logarithm" of concentration, and its exponential does not always exist, serving as a critical drawback both in theory and application. By contrast, the proposed model completely avoids this issue. Therefore, there is a tradeoff between the applicability and assumption of the models for the water quantity-quality dynamics. We report that there are major water quality indices that have positive correlations to discharge as discussed in **Section 4**.

## 3. Computational Aspects

Generating sample paths for the proposed model is essential for analyzing its probability density that has not been found analytically as well as event-wise concentration–discharge relationships. In this section, we discuss how to discretize our model with explanations of the difficulties faced while applying common numerical methods to the model and the reasons for using the method of Abi Jaber (2024)[60].

### 3.1 Difficulty in Simulating Square-root Process

The nominal square-root process $X=(X_t)_{t\geq 0}$ satisfies the SDE with an initial condition $X_0 > 0$:

$$dX_t = -R(X_t - f)dt + g\sqrt{RX_t}dW_t, \ t > 0 \tag{19}$$

with constants $R, f, g > 0$ and a one-dimensional standard Brownian motion $W=(W_t)_{t\in\mathbb{R}}$. The SDE (19) has a unique path-wise solution that is non-negative (Theorem 1.2.1 in Alfonsi (2015)[43]). The solution has distinct behaviors depending on values of $f, g$ (Chapter 2 in Alfonsi (2015)[43]); the solution $X$ is positive for a small diffusion such that $g^2 \leq 2f$, while it may become zero otherwise ($g^2 > 2f$). The latter case becomes a challenge in computationally sampling paths of the square-root process using a common numerical method such as the Euler–Maruyama method:

$$X_{t+\Delta t} = X_t - R(X_t - f)\Delta t + g\sqrt{R|X_t|}(W_{t+\Delta t} - W_t), \ t = k\Delta t \ (k=0,1,2,...) \tag{20}$$

with a sufficiently small increment in $\Delta t > 0$. The accuracy of the discretization (20) and its variants becomes arbitrarily poor as $g$ increases (Hefter and Herzwurm 2018; Mickel and Neuenkirch 2025)[44,45]. For our model, from Eq. (7), the accuracy of sample paths would become problematic if we formally have

$$\sigma^2 > 2(a+bQ_t)\rho(dR), \tag{21}$$

where "$\rho(\mathrm{d}R)$" here is regarded as a very small positive value considering it as a probability of taking values approximately at $R$ ($O(10^{-3})$ in **Section 4**) in computation. This condition can be violated as the computational resolution for approximating $\rho$ increases because $\sigma$ is a constant. Moreover, the non-negativity of the numerical solutions is not necessarily satisfied by naïve discretization Eq. (20) because $W_{t+\Delta t} - W_t$ can take any real value. These observations imply that common numerical methods are not applicable to the proposed model.

### 3.2 Dynamics-based Method

Abi Jaber (2024)[60] proposed a dynamics-based numerical method for generating sample paths of the SDE (19) such that numerical solutions can be computed explicitly as in the Euler–Maruyama method and their negativity is always satisfied. Moreover, the method applies to any values of $g \geq 0$ in a unified manner.

The discretization of Abi Jaber (2024)[60] for the SDE (19) is expressed as follows:

$$X_{t+\Delta t} = X_t + Rf\Delta t - RU_t + g\sqrt{R}Z_t, \ t = k\Delta t \ (k = 0, 1, 2, ...) \tag{22}$$

with

$$Z_t = \frac{1}{\omega_t}(U_t - \kappa_t), \tag{23}$$

$$U_t \sim IG\left(\kappa_t, \left(\frac{\kappa_t}{\omega_t}\right)^2\right), \tag{24}$$

$$\kappa_t = X_t \frac{1-e^{-R\Delta t}}{R} - f\left(\frac{1-e^{-R\Delta t}}{R} - \Delta t\right), \tag{25}$$

and

$$\omega_t = g\sqrt{R}\frac{1-e^{-R\Delta t}}{R}. \tag{26}$$

Herein, $IG(\mu, \eta)$ is the inverse Gaussian distribution with mean $\mu$ and probability density (PI is the circle ratio 3.1415…)

$$p_{IG}(x) = \sqrt{\frac{\eta}{2\mathrm{PI}x^3}}\exp\left(-\frac{\eta(x-\mu)^2}{2\mu^2 x}\right), \ x > 0. \tag{27}$$

Random variables following $IG(\mu, \eta)$ can be easily generated using an acceptance-rejection method. Each $U_t$ at $t = k\Delta t$ ($k = 0, 1, 2, ...$) is independent, and $IG(0, 0)$ generates a zero.

The form of the discretization approach can be derived from the consistency of moment-generating functions or Riccati equations associated with the original and discretized square-root processes, as in Section 1 by Abi Jaber (2024)[60]. For this study, the discretization method allows for any value of $\Delta t$, and numerical solutions always remain non-negative irrespective of the values of $g$ due to Theorem 1.3 by Abi Jaber (2024)[60].

### 3.3 Numerical Method for the Proposed Model

The discretization method presented above is adapted to our model. We assume that $Q$ to be sampled at each $t = k\Delta t$ while preserving its non-negativity, which is possible by numerically specifying the dynamics of $Q$ as a superposition process such as the superposition of Ornstein–Uhlenbeck processes (Fasen 2009; Yoshioka et al. 2023)[54,78]. Accordingly, we approximate the probability measure $\rho$ using an empirical measure under the assumption that $\rho$ admits a probability density:

$$\rho(\mathrm{d}R) \approx \frac{1}{N}\sum_{i=1}^{N} \delta(R - R_i). \tag{28}$$

Herein, $N \in \mathbb{N}$ is the degree of freedom of discretization and $R_i$ is the $\frac{2i-1}{2N}$ th quantile level of $\rho$: $\int_0^{R_i} \rho(\mathrm{d}R) = \frac{2i-1}{2N}$ (Yoshioka et al. 2023)[78]. We set the discretized memory process

$$M_t^{(N)} = \sum_{i=1}^{N} m_t^{(N,i)}, \ t \in \mathbb{R} \tag{29}$$

with

$$\mathrm{d}m_t^{(N,i)} = -R_i\left(m_t^{(N,i)} - \frac{1}{N}(a + bQ_t)\right)\mathrm{d}t + \sigma\sqrt{R_i m_t^{(N,i)}}\mathrm{d}W_t^i, \ i = 1,2,3,...,N. \tag{30}$$

Here, each $W_t^i$ is a one-dimensional standard Brownian motion that is mutually independent of the others and $Q$. The SDE (30) is a square-root process with a time-dependent drift coefficient. We discretize each SDE (30) directly using the method employed for Eq. (22) at each $t = k\Delta t$ ($k = 0,1,2,...$):

$$\kappa_t = m_t^{(N,i)}\frac{1 - e^{-R_i\Delta t}}{R} - \frac{1}{N}(a + bQ_t)\left(\frac{1 - e^{-R_i\Delta t}}{R_i} - \Delta t\right) \tag{31}$$

and

$$\omega_t = \sigma\sqrt{R_i}\frac{1 - e^{-R_i\Delta t}}{R_i}, \tag{32}$$

where all inverse Gaussian noises are sampled independently. The non-negativity of $m^{(N,i)}$ and $M^{(N)}$ is satisfied, as the proof of Theorem 1.3 by Abi Jaber (2024)[60] is applicable to our study owing to the non-negativity of $Q$. The computational performance of the discretization method is discussed in the next section using the identified parameter values.

### 3.4 Computational Tests

We examine the numerical method discussed in **Section 3.2** for three values of $f$ with $R = 1$ and $g = 10$. The cases considered in this paper include (A) $f = 100$ ($g^2 \le 2f$), (B) $f = 10$ ($g^2 > 2f$), and

(C) $f=0$. Case (C) is an extreme case where the square-root process becomes zero once it reaches state zero. Computational tests against the numerical method focus on the evaluation of long-term average and variance, along with path-wise properties, as all components are important for our application. For $t \geq 0$, by taking averages in the SDE (19), we obtain

$$\mathbb{E}[X_t] = X_0 e^{-t} + f(1-e^{-t}), \tag{33}$$

$$\mathbb{V}[X_t] = X_0 f^2 (e^{-t} - e^{-2t}) + \frac{fg^2}{2}(1-e^{-t})^2. \tag{34}$$

We numerically compute the average and variance for cases (A)–(C). Further, we numerically compute the average of the first hitting time $\tau$ for case (C). Here, we define

$$\tau = \inf\{t>0 | X_t = 0\}, \tag{35}$$

after which the process $X$ stays at zero. Abi Jaber (2024)[60] did not conduct such tests; particularly, the second test to determine the first hitting time aids in evaluating the path-wise performance of the numerical method as $\tau$ depends on the sample paths. The average of $\tau$ is available as a specific integral (Theorem 1.5 by Li et al. (2019)[79] with $q=0$, $x=X_0$, $b=f$, $c=g/2$, and $\theta=1$ in this literature):

$$\mathbb{E}[\tau] = \int_0^{+\infty} \frac{1-e^{-X_0 u}}{u(f+gu/2)} du. \tag{36}$$

We generate 1,000,000 sample paths for each computational test to determine the convergence of the discretization method. We examine different values of time increment $\Delta t$ to determine the convergence rate. Thus, the minimum digit ($\Delta t$) of $\tau$ for each sample path is $\Delta t$. The difference between the computed and theoretical average and variance are evaluated using least-squares errors during the computational period, which is fixed to $0 \leq t \leq 20$ for all computational experiments. The error between the theoretical and computed expectations of the first-hit time $\tau$ is measured by considering the absolute difference. The initial condition is fixed as $X_0 = 5$.

**Table 1** lists the errors in the computed average and variance of $X$ in each case. The average does not significantly change, but variance decreases with not so small $\Delta t$. **Table 2** shows the corresponding convergence rates. Sample paths for cases (A)–(C) are shown in **Fig. 1**. The path remains positive in case (A), sometimes close to be 0 in case (B), and exhibits the value 0 after a hitting time in case (C). The error in the average is attributed to the finiteness of the total number of sample paths because no errors existed in the average according to Proposition 1.4 by Abi Jaber (2024)[60]. The convergence rate is irregular and does not have clear trend, and that of the computed variance decreases with $f$ as diffusion is dominated. Moreover, the error does not decrease at the smallest $\Delta t = 0.0001$ examined for this study because the errors due to sample finiteness dominates for a sufficiently small $\Delta t$, as listed in the last row of **Table 1**. **Table 3** compares theoretical and computed averages of $\tau$ for case (C). Contrarily, doubling the total number of sample paths does not significantly improve the accuracy of the average of $\tau$ because the minimal digit to compute each $\tau$ is $\Delta t$. The error decreases with $\Delta t$, and the convergence rate is slightly lower than one. Using a smaller $\Delta t$ improves the accuracy of computed average of $\tau$ under the setting of this study as shown in **Table 3**.

The computed and theoretical average and variance suggest that the theoretical asymptotic values are reasonably captured at $\Delta t = 0.01$ as shown in **Fig. 2**. Considering a relatively large increment $\Delta t = 0.1$ yields an error of approximately 10 % in variance for cases (A) and (B). Therefore, the selection of $\Delta t$ is crucial for long-term numerical computations. The computational results obtained suggest that the total number of sample paths and computational resolution in time should be balanced depending on the objective of computation; in the present case, average and variance or first hitting time.

The results obtained suggest that the numerical method discussed in **Section 3.2** can deal with square-root processes for broad sizes of diffusion.

**Table 1** Errors in average and variance. "Path doubled" means that the total number of sample paths are doubled

| $\Delta t$ | Average | | | Variance | | |
|---|---|---|---|---|---|---|
| | (A) | (B) | (C) | (A) | (B) | (C) |
| 0.1 | 4.907.E-02 | 1.777.E-02 | 2.923.E-03 | 4.572.E+02 | 4.600.E+01 | 2.792.E+00 |
| 0.01 | 4.621.E-02 | 1.638.E-02 | 2.797.E-03 | 4.460.E+01 | 4.395.E+00 | 6.570.E-01 |
| 0.001 | 5.437.E-02 | 1.828.E-02 | 2.211.E-03 | 8.207.E+00 | 2.118.E+00 | 2.526.E-01 |
| 0.0001 | 4.562.E-02 | 1.155.E-02 | 2.151.E-03 | 5.976.E+00 | 1.418.E+00 | 3.330.E-01 |
| 0.0001 (Path doubled) | 3.056.E-02 | 8.352.E-03 | 1.260.E-03 | 4.774.E+00 | 9.948.E-01 | 2.311.E-01 |

**Table 2** Convergence rate corresponding to **Table 1**. The convergence rate at a given $\Delta t$ is calculated as $\log_{10}\left(\text{Error}_{10\Delta t} / \text{Error}_{\Delta t}\right)$, where $\text{Error}_{\Delta t}$ is the error with the time increment $\Delta t$.

| $\Delta t$ | Average | | | Variance | | |
|---|---|---|---|---|---|---|
| | (A) | (B) | (C) | (A) | (B) | (C) |
| 0.01 | 2.61.E-02 | 3.51.E-02 | 1.91.E-02 | 1.01.E+00 | 1.02.E+00 | 6.28.E-01 |
| 0.001 | -7.06.E-02 | -4.76.E-02 | 1.02.E-01 | 7.35.E-01 | 3.17.E-01 | 4.15.E-01 |
| 0.0001 | 7.62.E-02 | 2.00.E-01 | 1.18.E-02 | 1.38.E-01 | 1.74.E-01 | -1.20.E-01 |

**Table 3** Errors in the expectation of the first hitting time $\tau$. The exact value is 0.289273. The convergence rate at $\Delta t = 0.00005$ is computed as $\log_2\left(\text{Error}_{2\Delta t} / \text{Error}_{\Delta t}\right)$. The method of **Table 2** applies to the other cases.

| $\Delta t$ | Computed value | Error | Convergence rate |
|---|---|---|---|
| 0.1 | 9.203.E-01 | 6.310.E-01 | |
| 0.01 | 3.729.E-01 | 8.363.E-02 | 8.777.E-01 |
| 0.001 | 2.998.E-01 | 1.050.E-02 | 9.010.E-01 |
| 0.0001 | 2.906.E-01 | 1.353.E-03 | 8.901.E-01 |
| 0.00005 | 2.900.E-01 | 6.800.E-04 | 9.926.E-01 |
| 0.0001 (Path doubled) | 2.905.E-01 | 1.248.E-03 | |

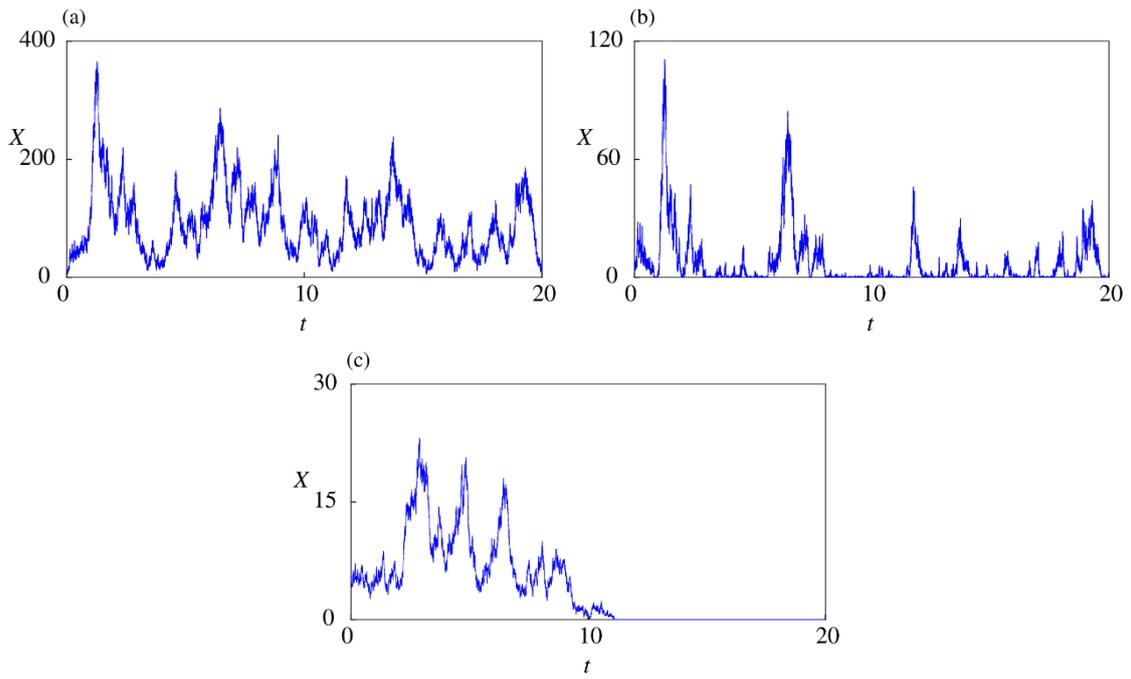

**Fig. 1** Sample paths of (a) case (A), (b) case (B), and (c) case (C)

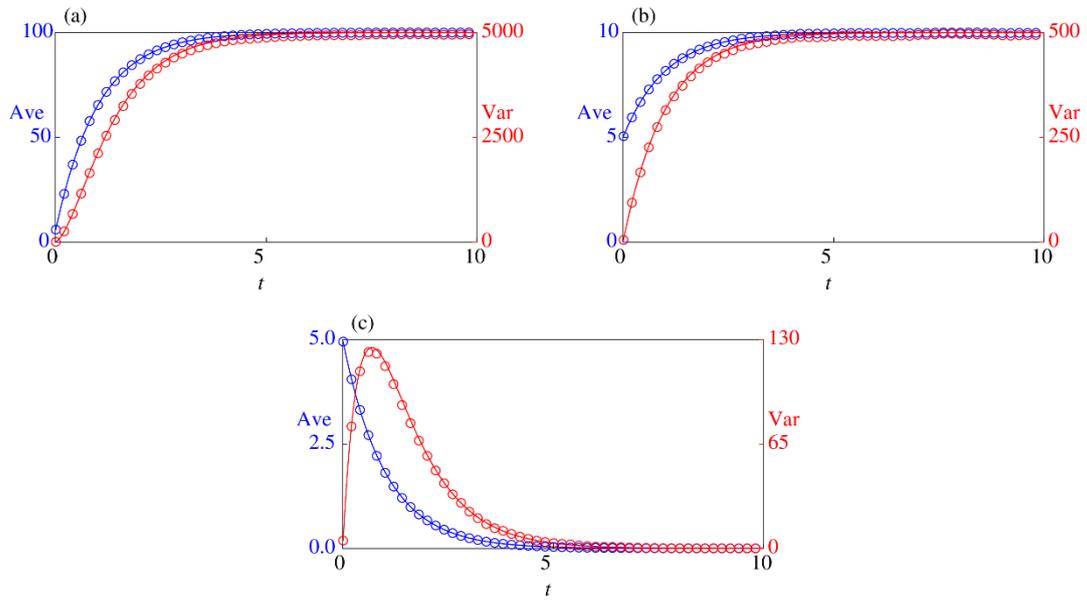

**Fig. 2** Computed and theoretical averages and variances for (a) case (A), (b) case (B), and (c) case (C). Line: Theoretical results. Circles: Computed results

## 4. Application

The proposed model is applied to water quality data sampled at a river in Japan.

### 4.1 Study Site and Variables

Water quality indices including TN, TP, and TOC were sampled almost weekly from the Hii River, flowing in the western part of Shimane Prefecture, Japan (**Fig. 3**), from 1991 to 2022 (Takeda 2023)[80]. The Hii River is 153 km long and has a watershed area of 2,540 km$^2$. The Otsu point has a rural catchment with an area of 911 km$^2$ mostly covered by forest (80.8 %). The remaining areas are covered by farmlands (10.0 %), urban areas (1.4 %), and others (7.8 %) (Takeda 2023)[80].

We use the time series data at the Otsu point to validate the proposed model because the public hourly discharge data were available from the water information system provided by (Ministry of Land, Infrastructure, Transport and Tourism 2025)[81]. A part of the water quality data was analyzed by Yoshioka and Yoshioka (2024a)[34] (data period 1991 to 2021). The data used in this study added samples during 2022 to the earlier samples and were made recently available. The total number of data points is 1,484 for both TN and TP, with the first and final sampling dates being August 20, 1991 and December 27, 2022, respectively. The total number of TOC data points is 820, with the first and final sampling dates being April 5, 2005 and December 27, 2022, respectively. The discharge data at Otsu corresponds to the period from 00:00 on January 1, 1991, to 23:00 on December 31, 2021. TN and TP were measured using spectrophotometry after potassium peroxodisulfate decomposition and molybdenum blue spectrophotometry, respectively. TOC was measured using the total organic carbon analyzer (SHIMADZU, TOC-Vcsn). The time-series data used in this study are plotted in **Figs. 4** and **5**.

Nitrogen and phosphorus are the major indices used to evaluate eutrophication and ecosystem states (Gillett et al. 2024; Liu et al. 2025; Luo et al. 2025)[73,82,83]. Monitoring their dynamics and concentration–discharge relationships aided in assessing anthropogenic impacts on freshwater ecosystems (Bieroza et al. 2024; Ehrhardt et al. 2021; Winter et al. 2020; Zhang et al. 2024)[84,85,86,87]. Fluvial carbon transport, evaluated using an index such as TOC, is essential for quantifying local and global carbon cycles and ecosystem functions. The dissolved organic carbon, a part of TOC, was positively correlated with discharge (Chen et al. 2021; Huntington and Shanley 2022; Drake et al. 2021; Drake et al. 2023; Goñi et al. 2023)[88,89,90,91,92].

In the Hii River, a system of cascading brackish lakes exists, referred to as Lake Shinji and Lake Nakaumi, which were registered as the Ramsar site. The Hii River, the largest river flowing into the connected lake system, has been important in this region (Hafeez and Inoue 2024; Ishitobi et al. 1998; Sugawara et al. 2017; Takeda et al. 2009)[93,94,95,96]; the Otsu point is close to the mouth of the connected lake system (**Fig. 3**), which is why we have chosen this observation point.

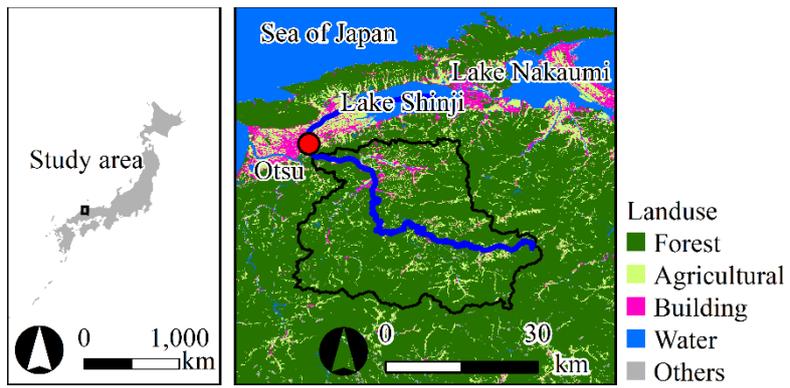

**Fig. 3** Hii River and the catchment of Otsu point

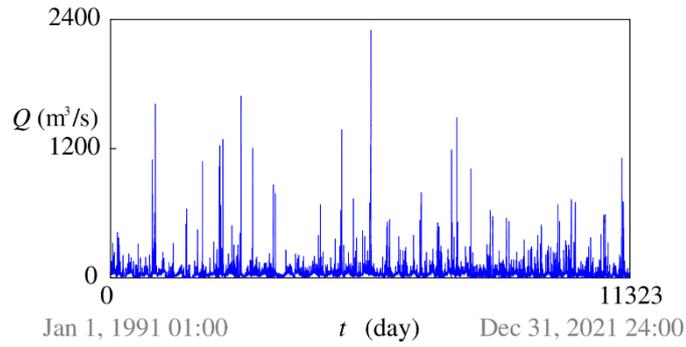

**Fig. 4** Time series of discharge at Otsu point

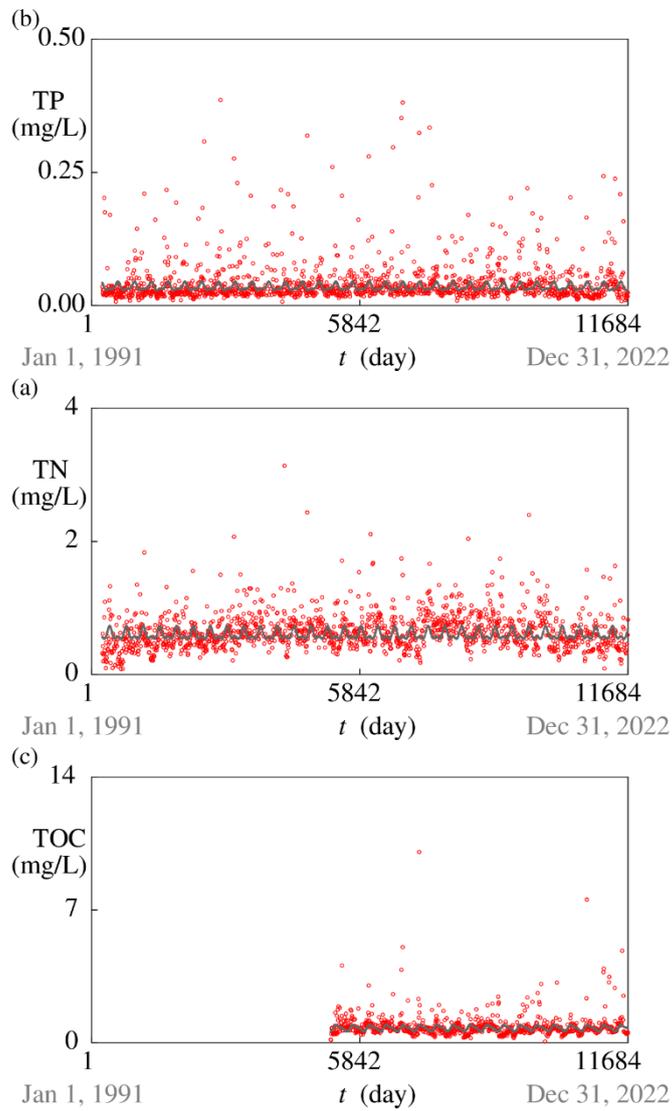

**Fig. 5** Time series of (a) TN, (b) TP, and (c) TOC: empirical data (red) and fitted seasonal component (grey). TN: 5.984 (mg/L) and TP: 0.818 (mg/L) were recorded on June 29, 1993 (day 911)

## 4.2 Model Identification

### 4.2.1 Discharge

The discharge model is specified as a superposition of the Ornstein–Uhlenbeck processes, as has been applied to other points in the Hii River and other rivers (Yoshioka and Yoshioka 2024b; Yoshioka and Yoshioka 2025)[37,38]:

$$Q_t = \int_{s=-\infty}^{s=t} \int_{z=0}^{z=+\infty} \int_{r=0}^{r=+\infty} z e^{-r(t-s)} \Xi(\mathrm{d}s, \mathrm{d}z, \mathrm{d}r), \ t \in \mathbb{R}, \quad (37)$$

where $z > 0$ is the jump size of discharge, $r > 0$ is the recession rate of each flood pulse, and $\Xi(\mathrm{d}s, \mathrm{d}z, \mathrm{d}r)$ is a Poisson random measure of the three-dimensional phase-space $\mathbb{R} \times (0, +\infty) \times (0, +\infty)$ of time, jump size, and recession rate. The compensator associated with $\Xi(\mathrm{d}s, \mathrm{d}z, \mathrm{d}r)$ is $\mathrm{d}s \nu(\mathrm{d}z) \pi(\mathrm{d}r)$, with the jump intensity measure of the form $\nu(\mathrm{d}z) = a_1 e^{-a_2 z} z^{-(1+a_3)} \mathrm{d}z$, where $a_1, a_2 > 0$ and $a_3 < 1$. The model (37) assumes that the slow (small $r$, groundwater) to fast (large $r > 0$, surface water) runoffs contribute to discharge.

The probability measure $\pi(\mathrm{d}r)$ of the recession rate $r$ is assumed to be gamma type

$$\pi(\mathrm{d}r) = \frac{1}{\Gamma(\alpha)\beta^\alpha} r^{\alpha-1} e^{-r/\beta} \mathrm{d}r, \ r > 0 \quad (38)$$

with $\alpha > 1$ and $\beta > 0$ that can cover both slow and fast runoffs and moreover gives the autocorrelation of the discharge explicitly as $A_Q(h) = (1+\beta h)^{-(\alpha-1)}$ ($h \geq 0$) with the following $\theta$:

$$\theta(\mathrm{d}r) = \frac{r^{-1} \pi(\mathrm{d}r)}{\int_0^{+\infty} r^{-1} \pi(\mathrm{d}r)}, \ r > 0. \quad (39)$$

These modeling specifications have been successfully applied to the existing discharge data by employing the identification method explained below (Yoshioka and Yoshioka 2024b; Yoshioka and Yoshioka 2025)[37,38]. Firstly, the parameters $\alpha$ and $\beta$ are identified by minimizing the least-squares error between empirical and theoretical autocorrelations. Secondly, the remaining model parameters $a_1, a_2, a_3$ are determined using the moment-matching method to minimize the sum of relative errors:

$$\left(\frac{\bar{Q}_e - \bar{Q}_m}{\bar{Q}_e}\right)^2 + \left(\frac{\bar{V}_e - \bar{V}_m}{\bar{V}_e}\right)^2 + \left(\frac{\bar{S}_e - \bar{S}_m}{\bar{S}_e}\right)^2 + \left(\frac{\bar{K}_e - \bar{K}_m}{\bar{K}_e}\right)^2, \quad (40)$$

where $\bar{S}$ and $\bar{K}$ represent the skewness and kurtosis, respectively, both of which exist owing to the assumptions $\nu$. The subscripts "e" and "m" represent empirical and modelled quantities, respectively. Moment matching is employed in this study because discharge data have high skewness and kurtosis, larger than $O(10^0)$ due to spikes corresponding to flood events.

### 4.2.2 Water Quality

Having specified the model of discharge $Q$, we identify the model of concentration $C$. The identification method is based on the approach proposed by Yoshioka and Yoshioka (2025)[38] with modifications such that both short and long memories are dealt with; the original method focused only on long memories.

The one-year period $T$ is set as 365.25 days. We assume a sinusoidal model as the periodicity:

$$S_t = e^{A_0 + A_1 \sin\left(\frac{2\mathrm{PI}t}{T} + B_1\right) + A_2 \sin\left(\frac{4\mathrm{PI}t}{T} + B_2\right)} \tag{41}$$

with time zero at the beginning of January of each year. The five real parameters $A_0, A_1, A_2, B_1, B_2$ are identified by least-squares fitting between the empirical and theoretical $\ln(C_t)$. Subsequently, $M$ is obtained using $M_t = C_t S_t^{-1}$, with the right-hand side computed using the identified parameter values.

Further, we identify the parameters $a, b, \sigma$ and the probability measure $\rho$. By Eq. (11), the autocorrelation $A_M(h)$ is rewritten as:

$$A_M(h) = \frac{\int_0^{+\infty} e^{-Rh} \rho(\mathrm{d}R) + w(I(h) + J(h))}{1 + wI(0)}, \quad h \geq 0 \tag{42}$$

with $w = \frac{2\mu^2 \bar{V}_e}{\sigma^2 \bar{M}_e} > 0$, $J(0) = 0$, and the notation $\bar{M}_e = \mathbb{E}[M_t]_e$. Assuming a gamma-type probability measure $\rho$ of Eq (38), with shape parameter $\varsigma > 0$ and scale parameter $\xi > 0$ as the simplest model as in the case for discharge, we temporally fit their values by setting $w = 0$. Hence, $A_M(h) = (1 + \xi h)^{-\varsigma}$ against the least-squares fitting for empirical autocorrelation is derived. With this fitting procedure, we can not only estimate the tail of the autocorrelation but also obtain the nominal values of $\varsigma = \varsigma_n$ and $\xi = \xi_n$ used in the sequel.

We account for correlation between discharge and concentration in the parameter identification as follows. We fix the parameter value $\xi = \xi_n$, examine different values of $\varsigma$ around $\varsigma_n$, compute the right-hand side terms of Eq. (42), and determine $w$ by least-squares fitting between the theoretical Eq. (42) and empirical correlations. Then, we determine $\sigma$ as follows:

$$\sigma = \sqrt{\frac{2\mathbb{V}[M_t]_e}{\bar{M}_e (1 + wI(0))}}, \tag{43}$$

which is obtained from the relationship

$$\mathbb{V}[M_t] = \frac{1}{2} \sigma^2 \mathbb{E}[M_t] (1 + wI(0)), \tag{44}$$

where each cumulant can be replaced with an empirical cumulant. Subsequently, we obtain $b$ from

$$b = \sigma \sqrt{\frac{w \bar{M}_e}{2 \bar{V}_e}}, \tag{45}$$

where the right-hand side term is known at this stage. Finally, we obtain $a$ from Eq. (9):

$$a = b\bar{Q}_e - \bar{M}_e. \tag{46}$$

The procedure is repeated so that the following difference between empirical and theoretical covariance is minimized:

$$\left| \lambda(0) \right|_e - \left| \lambda(0) \right|_m \right|. \tag{47}$$

The empirical covariance $\left.\lambda(0)\right|_e$ is computed using the measured concentration and daily averaged discharge on each sampling date. The manual component in the procedure examined is to determine values of $\varsigma$ around $\varsigma = \varsigma_n$. The calculation is conducted with increments of 0.0005 for TN, which is suggested to have a long memory. However, preliminary investigations suggested that TP and TOC exhibit autocorrelations that were almost exponential. In these cases, we use $\rho(dR) = \delta(R - \xi)$ instead of the gamma distribution and adapt the identification procedure against $\xi$. The increment for exploring the best $\xi$ is set as 0.001. The parameter $\varsigma$ does not appear in the model of TP and TOC.

Each integral in the variance, autocorrelation, and covariance is discretized by the quantile integral explained in **Section 3.3** with resolution $N = 2,048$. This resolution is used in the sample path generation of our model. Probability densities of discharge, TN, TP, and TOC are unavailable and are therefore numerically computed by sampling a 1,000-year sample path after a 1,000-year period for burn-in. The time increment in the Monte Carlo simulation is 0.01 (day). The model (37) of the discharge $Q$ is discretized using the Euler–Maruyama discretization with the same $N$, and the numerical method in **Section 3.3** is applied to the discretization of memory process $M$.

### 4.3 Fitting Results

The identified parameter values are summarized in **Tables 4** and **5** and a comparison of the theoretical and empirical statistics is presented in **Tables 6** and **7**. The theoretical and empirical probability densities of discharge, TN, TP, and TOC at the Otsu point are shown in **Fig. 6**, while **Fig. 7** shows the same results on an ordinary logarithmic scale. The current computational resolution of the Monte Carlo simulation, along with our numerical method is reasonable, as discussed in **Section A3**. The computed statistics reasonably align with the theoretical results despite the existence of errors in the discretization of the time derivative and probability measures $\rho, \theta$ in the proposed model, along with a finite number of samples. The fitting in the logarithmic scale captures the empirical tails for all the variables.

The discharge exhibits a unimodal density, which is captured by the identified model. Moreover, the tail of the empirical density is well fitted by the identified model, and the empirical cumulants are reproduced with a relative error of less than 0.003. The concentrations of TN, TP, and TOC at Otsu exhibit unimodal histograms, which are consistent with the literature on different surface water systems (Xie et al. 2022, Duan et al. 2025, Wan et al. 2024, Perri and Porporato 2022 [97,98,99,100]; the last study dealt with the sum of nitrate and nitrite). The theoretical and empirical autocorrelations of discharge and memory processes of TN, TP, and TOC are shown in **Fig. 8**, suggesting a reasonable agreement, with slight discrepancies for small time lags.

**Table 4** Identified parameter values for the model of discharge

| Parameter | Value |
|---|---|
| $\alpha$ (-) | 1.752.E+00 |
| $\beta$ (1/day) | 1.608.E+00 |
| $a_1$ ( $m^{3a_3}/s^{a_3}/day$ ) | 2.985.E+00 |
| $a_2$ ($s/m^3$) | 1.510.E-03 |
| $a_3$ (-) | 7.998.E-01 |

**Table 5** Identified parameter values for the models of TN, TP, and TOC

| | TN | TP | TOC |
|---|---|---|---|
| $e^{A_0}$ (mg/L) | 5.975.E-01 | 3.432.E-02 | 7.598.E-01 |
| $A_1$ (-) | 1.200.E-01 | -1.783.E-01 | -1.942.E-01 |
| $A_2$ (-) | 7.878.E-02 | 9.344.E-02 | 7.817.E-02 |
| $B_1$ (-) | 1.071.E+00 | 1.362.E+00 | 4.759.E-01 |
| $B_2$ (-) | 6.825.E-01 | 5.727.E-01 | 1.468.E+00 |
| $\varsigma_n$ (-) | 6.676.E-01 | | |
| $\xi_n$ (1/day) | 3.482.E-01 | 1.513.E-01 | 1.011.E-01 |
| $\varsigma$ (-) | 5.355E-01 | | |
| $\xi$ (1/day) | 3.482.E-01 | 2.940.E-01 | 2.380.E-01 |
| $\sigma$ (-) | 5.412.E-01 | 1.132.E+00 | 4.502.E-01 |
| $a$ (-) | 3.844.E-01 | 7.064.E-03 | 1.966.E-01 |
| $b$ ($m^3/s$) | 1.684.E-02 | 2.986.E-02 | 2.220.E-02 |

**Table 6** Comparison of theoretical and empirical statistics for the model of discharge $Q$

| | Empirical | Theoretical | Relative error |
|---|---|---|---|
| Average ($m^3/s$) | 4.1556.E+01 | 4.1560.E+01 | 9.2192.E-05 |
| Variance ($m^6/s^2$) | 2.7549.E+03 | 2.7548.E+03 | 4.4706.E-05 |
| Skewness (-) | 1.0073.E+01 | 1.0096.E+01 | 2.2398.E-03 |
| Kurtosis (-) | 2.1042.E+02 | 2.1019.E+02 | 1.0719.E-03 |

**Table 7** Comparison of theoretical and empirical statistics for the memory processes $M$ of TN, TP, and TOC

| | TN | | TP | | TOC | |
|---|---|---|---|---|---|---|
| Average (-) | 1.084.E-01 | | 1.248.E+00 | | 1.119.E+00 | |
| Variance (-) | 2.676.E-01 | | 1.659.E+00 | | 5.432.E-01 | |
| Covariance $\lambda(0)$ | Empirical | Theoretical | Empirical | Theoretical | Empirical | Theoretical |
| ($mg \cdot m^3/L/s$) | 9.778.E+00 | 9.766.E+00 | 2.879.E+01 | 2.879.E+01 | 1.937.E+01 | 1.936.E+01 |

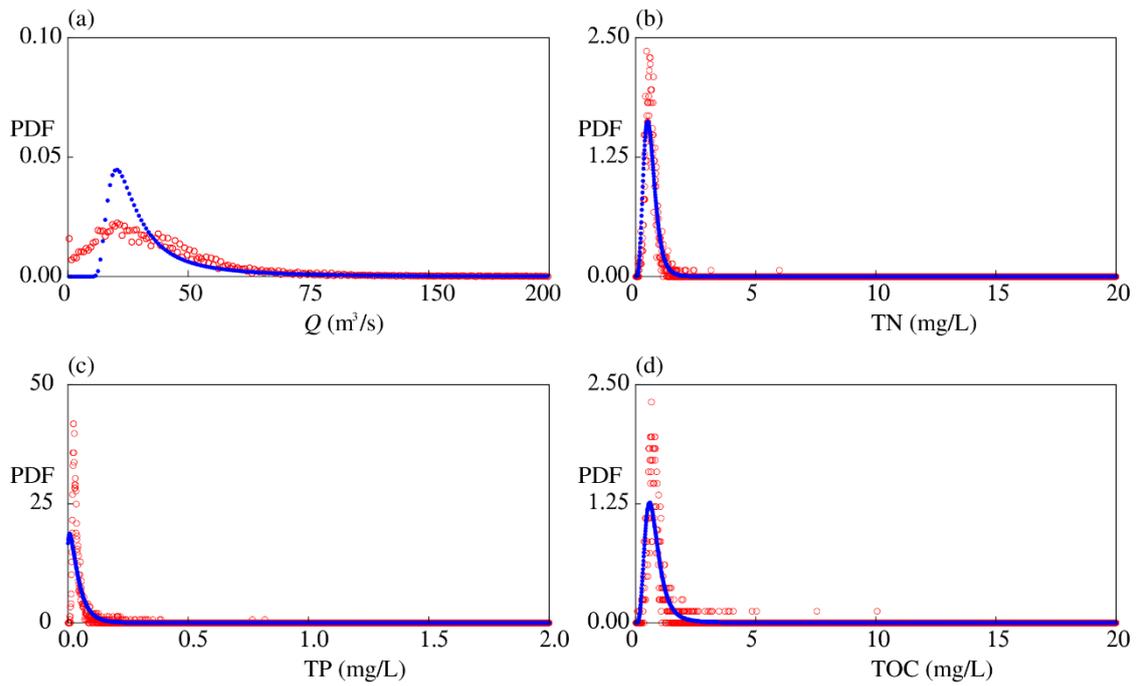

**Fig. 6** Comparison of theoretical and empirical probability density functions (PDFs) of (a) discharge $Q$ (m$^3$/s), (b) TN, (c) TP, and (d) TOC, showing empirical (red) and computed results (blue)

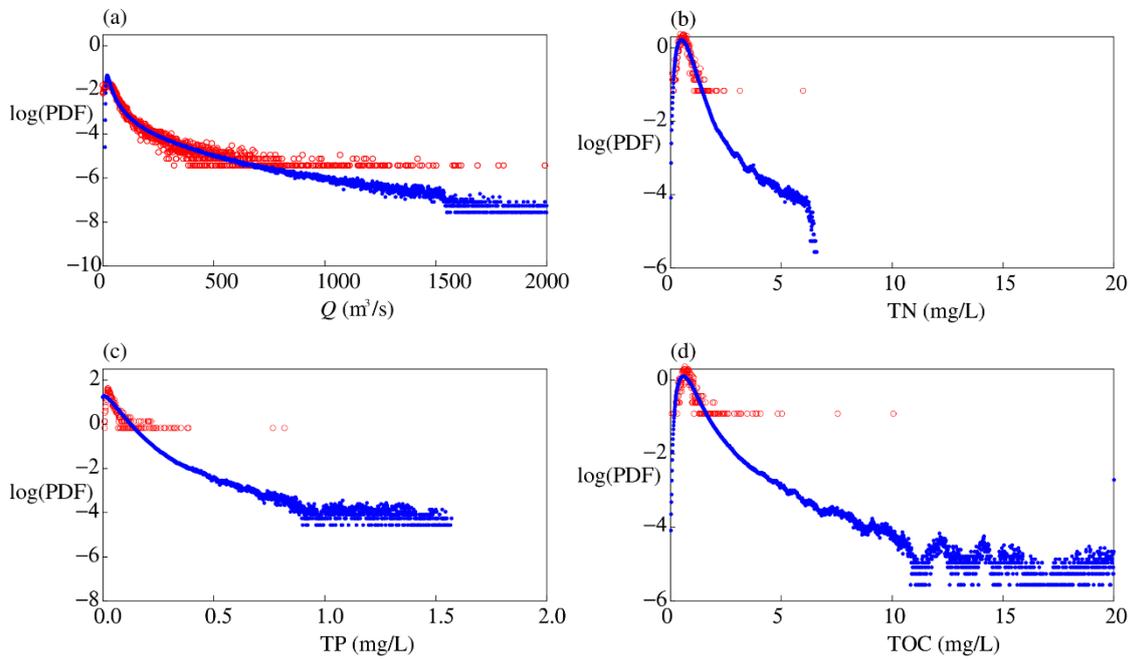

**Fig. 7** Same results as those shown in **Fig. 6** but on the ordinary logarithmic scale

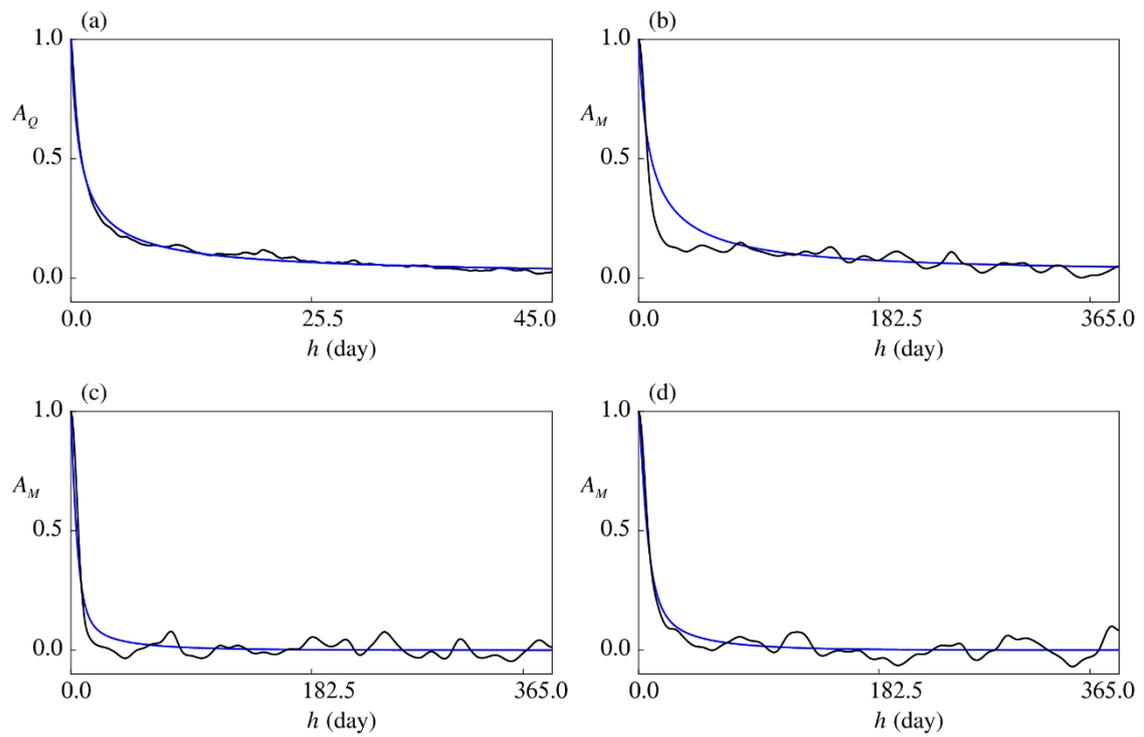

**Fig. 8** Theoretical and empirical autocorrelations of (a) discharge, (b) TN, (c) TP, and (d) TOC, showing empirical (black) and theoretical results (blue)

## 4.4 Concentration–discharge Relationship

First, we analyze the mutual covariance $\lambda$ between discharge and TN, TP, and TOC at the Otsu point, based on **Proposition 3**. Each $\lambda(h)$ for $h \geq 0$ with normalization $\lambda(0) = 1$ is plotted in **Fig. 9**. The negative branch for $h < 0$ is not plotted here because it has been proven to be monotonic in **Section 3**. The computed mutual covariances are maximized at the lag of 1.07, 1.20, and 1.42 days for TN, TP, and TOC, respectively, suggesting that the peak of these water quality indices appears approximately one-day later than that of discharge. The daily scale information is not directly found from the time series data in this study because the water quality is measured only weekly, demonstrating the advantage of using the proposed model to statistically evaluate the hysteresis in the concentration–discharge relationship.

The sample paths of discharge and TN, TP, and TOC are shown in **Fig. 10**. The paths of TN and TOC are similar to that of the square-process with a sufficiently small noise intensity (case (A) in **Section 3**), while that of TP is similar to that with a large noise intensity (case (B) in **Section 3**). This qualitative difference is attributed to the difference in the parameter $\sigma$; TP has a larger $\sigma$ than those of TN and TOC. Examples of the concentration–discharge relationship of TN are given in **Fig. 11**, showing that manually identifying rising and falling limbs to apply the classical hysteresis indices is not simple. By contrast, the approach based on mutual covariance completely avoids this technical issue.

Finally, we investigate the influence of model misspecification on mutual covariance, particularly its peak location. This analysis has been motivated by the fact that completely fitting a mathematical model to real data is difficult in general . We focus on the TN model owing to a larger number of parameters than TP and TOC because of the use of non-atomic $\rho$. We consider scenarios where the identified scale parameter $\xi$ is misspecified as $\xi(1+\Delta)$, with some $\Delta > -1$. The difference between the models with and without misspecification is evaluated by the Kullback–Leibler (KL) divergence as the simplest measure of the difference between the two probability measures (Kullback and Leibler 1951)[101]. After elementary calculations, the KL divergence as a function of $\Delta$ is obtained as follows:

$$\mathrm{KL}(\Delta) = \varsigma(\Delta - \ln(1+\Delta)), \tag{48}$$

which is convex and had a strict minimum value of zero in the $\Delta = 0$. The model with the identified $\xi$ is regarded as the baseline model, and that with $\xi(1+\Delta)$ is considered the alternative model. Theoretically, other divergence models such as the *f*-divergence (Sason and Verdú, 2016)[102] can be used as preferred; however, such model may introduce parameters that would require tuning.

**Table 8** lists the peak location, corresponding KL divergence, and the covariance $\lambda(0)$ for different values of $\Delta$. The mutual covariance $\lambda = \lambda(h)$ as a function of lag $h \geq 0$ (day) for different values of $\Delta$ is shown in **Fig. 12**. The peak location of the mutual covariance moves to the left as $\Delta$ increases, showing that an overestimation of the scale parameter $\xi$ results in a shorter delay between the peaks of discharge and TN on average, and vice versa. In addition, the covariance increases as $\xi$ decreased. Mathematically, the increase of net reversion with Δ causes the increase in covariance, according to the SDE (4). From **Table 8**, introducing 20 % errors ( $\Delta = \pm 0.2$ ) results in KL divergence of approximately 0.1, with at most 10 % difference in covariance between discharge and TN. In this way, we can determine the relations between size of misspecification and their influence on statistics. The high tractability of the proposed model allows for the evaluation of covariance and mutual covariance corresponding to long and short perspectives of the concentration–discharge relationships.

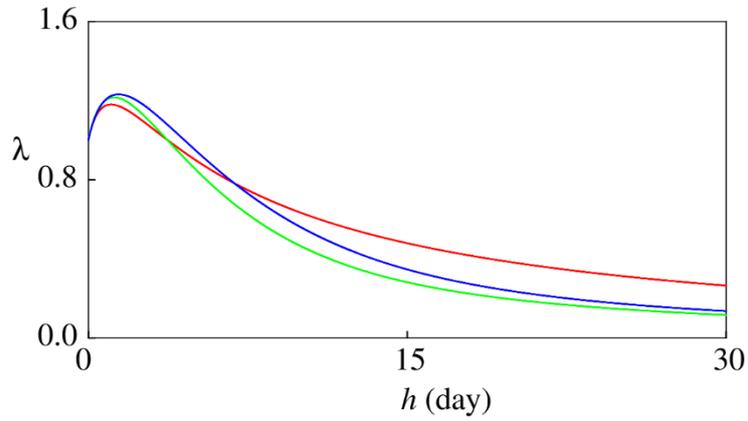

**Fig. 9** Mutual covariance $\lambda = \lambda(h)$ as a function of lag $h \geq 0$, between discharge and TN (red), discharge and TP (green), and discharge and TOC (blue)

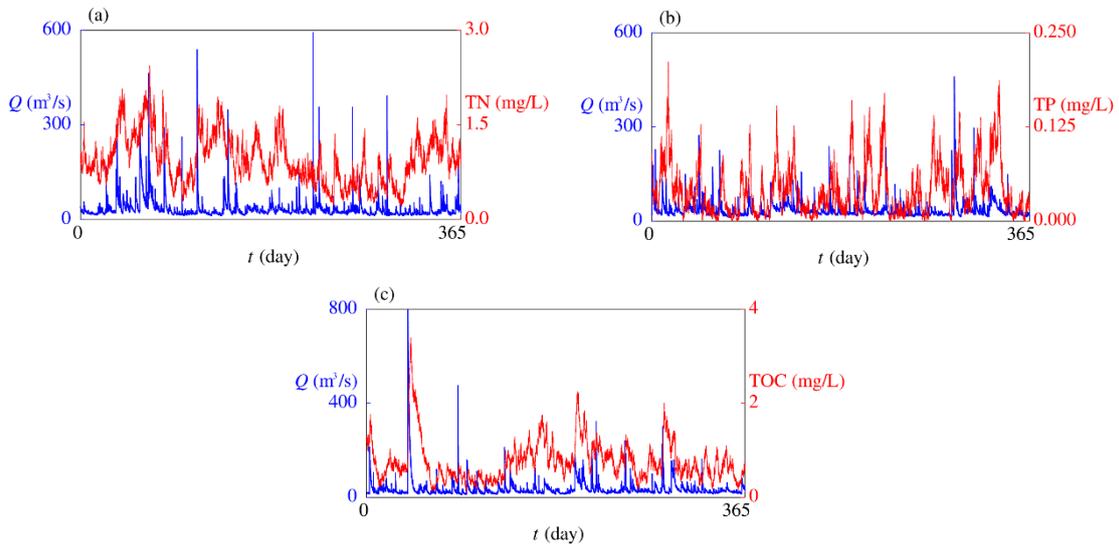

**Fig. 10** Sample paths: (a) discharge and TN, (b) discharge and TP, and (c) discharge and TOC

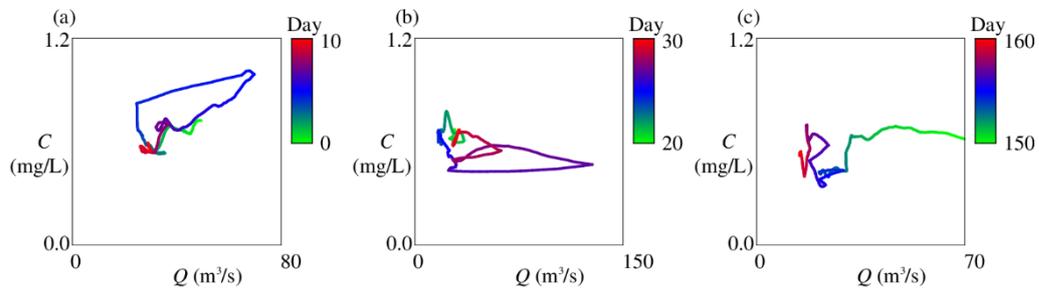

**Fig. 11** Concentration–discharge relationship of TN at selected time intervals for the time series presented in **Fig. 10(a)**. Daily averaged discharge and TN concentration are used for visualization

**Table 8** Peak location and corresponding KL divergence for different values of $\Delta$

| $\Delta$ | Peak location (day) | KL($\Delta$) | $\lambda(0)$ |
|---|---|---|---|
| -0.5 | 1.74 | 1.034.E-01 | 6.794.E+00 |
| -0.4 | 1.53 | 5.935.E-02 | 7.508.E+00 |
| -0.3 | 1.37 | 3.035.E-02 | 8.149.E+00 |
| -0.2 | 1.25 | 1.239.E-02 | 8.733.E+00 |
| -0.1 | 1.15 | 2.871.E-03 | 9.269.E+00 |
| 0 | 1.07 | 0.000.E+00 | 9.766.E+00 |
| 0.1 | 1.00 | 2.511.E-03 | 1.023.E+01 |
| 0.2 | 0.94 | 9.467.E-03 | 1.066.E+01 |
| 0.3 | 0.89 | 2.015.E-02 | 1.107.E+01 |
| 0.4 | 0.85 | 3.402.E-02 | 1.145.E+01 |
| 0.5 | 0.81 | 5.062.E-02 | 1.182.E+01 |

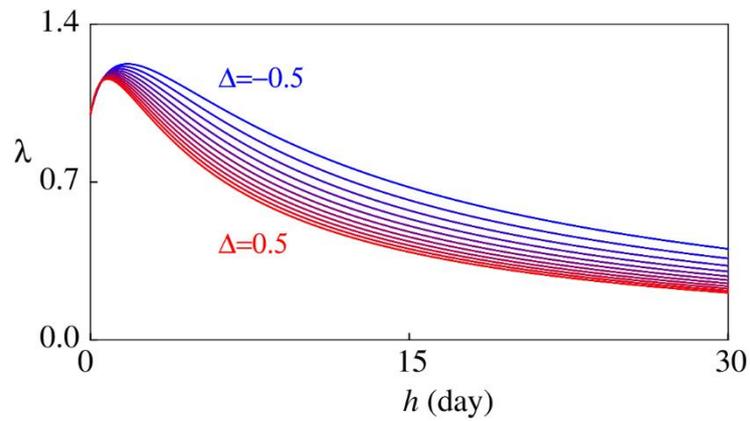

**Fig. 12** Mutual covariance $\lambda = \lambda(h)$ as a function of lag $h \geq 0$ (day) for different values of $\Delta$ from $-0.5$ (blue, top) to $0.5$ (red, bottom) at uniform intervals of $0.1$. Normalization is used such that $\lambda(h) = 1$ in the absence of lag ($h = 0$)

## 5. Conclusions

This study introduced a non-Markovian mathematical model designed to address both short- and long-term concentration–discharge relationships, along with its statistics such as average, variance, autocorrelation, and mutual covariance. To this end, a numerical method for simulating the sample paths of this model was introduced, based on a recently developed discretization technique for generic square-root processes. The numerical method successfully managed the multi-scale nature of the proposed model over time, while keeping the non-negativity of the numerical solutions. The case study suggested that the proposed model could be applied to water quality indices exhibiting both short and long memories. This study quantified the concentration–discharge relationship using mutual covariance, indicating an approximate one-day delay in the peaks of TN, TP, and TOC relative to discharge. Additionally, a sensitivity analysis was performed to evaluate the effect of model uncertainty on mutual covariance.

The results obtained in this study suggested that our mathematical model can be applied to other river environments and water quality indices if the identification method is successful and sufficient data are collected. Theoretically, the model and its mutual covariance with proper modifications can be applied to analyze the hysteresis between surface water and groundwater dynamics, along with the water stage-discharge dynamics in rivers. The exact statistics obtained in this paper are novel contributions by themselves; moreover, they can be used for verifying numerical methods for superposition processes.

Furthermore, the versatility of the numerical method proposed by Abi Jaber (2024)[60] was demonstrated in this study. We anticipate that this method will serve as a breakthrough for problems where simulating square root and related processes with large diffusion coefficients is critical. Future applications of this method include analyzing fish migration, where the dynamics of migrating fish populations are influenced by water quantity and temperature variations. Additionally, an intriguing avenue for future research involves optimizing river environments through intensive human interventions, with the proposed model serving as a tool in the system dynamics to be controlled.

# Appendix

## A1. Proofs

*Proof of Proposition 1*

First, to prove Eq. (9), we take the expectation of Eq. (3) to obtain

$$\begin{aligned}
\mathbb{E}[M_t] &= \mathbb{E}\left[\int_0^{+\infty} m_t(\mathrm{d}R)\right] \\
&= \mathbb{E}\left[\int_{R=0}^{R=+\infty}\int_{s=-\infty}^{s=t} e^{-R(t-s)} R(a+bQ_s)\rho(\mathrm{d}R)\mathrm{d}s + \int_{R=0}^{R=+\infty}\int_{s=-\infty}^{s=t} e^{-R(t-s)}\sigma\sqrt{R} B(\mathrm{d}R,\mathrm{d}s)\right] \\
&= \int_{R=0}^{R=+\infty}\int_{s=-\infty}^{s=t} e^{-R(t-s)} R(a+b\mathbb{E}[Q_s])\rho(\mathrm{d}R)\mathrm{d}s \\
&= (a+b\bar{Q})\int_{R=0}^{R=+\infty}\int_{s=-\infty}^{s=t} e^{-R(t-s)} R\rho(\mathrm{d}R)\mathrm{d}s \\
&= (a+b\bar{Q})\int_{R=0}^{R=+\infty} R\rho(\mathrm{d}R)\int_{s=-\infty}^{s=t} e^{-R(t-s)}\mathrm{d}s \\
&= a+b\bar{Q}
\end{aligned} \quad , t\in\mathbb{R} \quad (A1)$$

because $B$ is Gaussian with mean 0.

To prove Eq. (10), we use the definition of variance and Eq. (9) to obtain

$$\begin{aligned}
\mathbb{V}[M_t] &= \mathbb{E}\left[(M_t - \mathbb{E}[M_t])^2\right] \\
&= \mathbb{E}\left[\left(\int_{R=0}^{R=+\infty}\int_{s=-\infty}^{s=t} e^{-R(t-s)}\left\{bR(Q_s-\bar{Q})\rho(\mathrm{d}s) + \sigma\sqrt{R}B(\mathrm{d}R,\mathrm{d}s)\right\}\right)^2\right].
\end{aligned} \quad (A2)$$

Applying isometry to Eq. (A2) along with the covariance structure (5), we obtain

$$\begin{aligned}
\mathbb{V}[M_t] &= b^2 \mathbb{E}\left[\int_{R=0}^{R=+\infty}\int_{s=-\infty}^{s=t}\int_{P=0}^{P=+\infty}\int_{u=-\infty}^{u=t} e^{-R(t-s)-P(t-u)} PR(Q_s-\bar{Q})(Q_u-\bar{Q})\rho(\mathrm{d}R)\rho(\mathrm{d}P)\mathrm{d}s\mathrm{d}u\right] \\
&\quad + \sigma^2 \mathbb{E}\left[\int_{R=0}^{R=+\infty}\int_{s=-\infty}^{s=t} Re^{-2R(t-s)} m_s(\mathrm{d}R)\mathrm{d}s\right].
\end{aligned} \quad (A3)$$

Each expectation in Eq. (A3) can be rewritten as:

$$\begin{aligned}
\mathbb{E}\left[\int_{R=0}^{R=+\infty}\int_{s=-\infty}^{s=t} Re^{-2R(t-s)} m_s(\mathrm{d}R)\mathrm{d}s\right] &= \mathbb{E}\left[\int_{R=0}^{R=+\infty}\int_{s=-\infty}^{s=t} Re^{-2R(t-s)} m_t(\mathrm{d}R)\mathrm{d}s\right] \\
&= \mathbb{E}\left[\int_{R=0}^{R=+\infty}\left(\int_{s=-\infty}^{s=t} Re^{-2R(t-s)}\mathrm{d}s\right) m_t(\mathrm{d}R)\right] \\
&= \frac{1}{2}\mathbb{E}\left[\int_0^{+\infty} m_t(\mathrm{d}R)\right] \\
&= \frac{1}{2}\mathbb{E}[M_t]
\end{aligned} \quad (A4)$$

by the stationarity of $m$, and

$$\mathbb{E}\left[\int_{R=0}^{R=+\infty}\int_{s=-\infty}^{s=t}\int_{P=0}^{P=+\infty}\int_{u=-\infty}^{u=t}e^{-R(t-s)-P(t-u)}PR(Q_s-\bar{Q})(Q_u-\bar{Q})\rho(\mathrm{d}R)\rho(\mathrm{d}P)\mathrm{d}s\mathrm{d}u\right]$$

$$=\int_{R=0}^{R=+\infty}\int_{s=-\infty}^{s=t}\int_{P=0}^{P=+\infty}\int_{u=-\infty}^{u=t}PRe^{-R(t-s)-P(t-u)}\mathbb{E}\left[(Q_s-\bar{Q})(Q_u-\bar{Q})\right]\rho(\mathrm{d}R)\rho(\mathrm{d}P)\mathrm{d}s\mathrm{d}u$$

$$=\int_{R=0}^{R=+\infty}\int_{s=-\infty}^{s=t}\int_{P=0}^{P=+\infty}\int_{u=-\infty}^{u=t}PRe^{-R(t-s)-P(t-u)}\bar{V}A_Q(|s-u|)\rho(\mathrm{d}R)\rho(\mathrm{d}P)\mathrm{d}s\mathrm{d}u$$

$$=\bar{V}\int_{R=0}^{R=+\infty}\int_{s=-\infty}^{s=t}\int_{P=0}^{P=+\infty}\int_{u=-\infty}^{u=t}\int_{r=0}^{r=+\infty}PR\left(e^{-R(t-s)-P(t-u)-r|s-u|}\theta(\mathrm{d}r)\right)\rho(\mathrm{d}R)\rho(\mathrm{d}P)\mathrm{d}s\mathrm{d}u \qquad (A5)$$

$$=\bar{V}\int_{R=0}^{R=+\infty}\int_{P=0}^{P=+\infty}\int_{r=0}^{r=+\infty}PR\left(\int_{s=-\infty}^{s=t}\int_{u=-\infty}^{u=t}e^{-R(t-s)-P(t-u)-r|s-u|}\mathrm{d}s\mathrm{d}u\right)\theta(\mathrm{d}r)\rho(\mathrm{d}R)\rho(\mathrm{d}P)$$

$$=\bar{V}\int_{R=0}^{R=+\infty}\int_{P=0}^{P=+\infty}\int_{r=0}^{r=+\infty}\frac{PR}{P+R}\left(\frac{1}{P+r}+\frac{1}{R+r}\right)\theta(\mathrm{d}r)\rho(\mathrm{d}R)\rho(\mathrm{d}P)$$

Here, we used $\mathbb{E}\left[(Q_s-\bar{Q})(Q_u-\bar{Q})\right]=\bar{V}A_Q(|s-u|)$ and the following elementary calculations:

$$\int_{-\infty}^{t}e^{-2R(t-s)}R\mathrm{d}s=\frac{1}{2}, \qquad (A6)$$

$$\int_{s=-\infty}^{s=t}\int_{u=-\infty}^{u=t}e^{-R(t-s)-P(t-u)-r|t-s|}\mathrm{d}s\mathrm{d}u=\int_{s=-\infty}^{s=0}\int_{u=-\infty}^{u=0}e^{Rs+Pu-r|s-u|}\mathrm{d}s\mathrm{d}u=\frac{1}{P+R}\left(\frac{1}{P+r}+\frac{1}{R+r}\right). \qquad (A7)$$

By substituting Eqs. (A4)-(A5) into Eq. (A3), we complete the proof.

□

*Proof of Proposition 2*

Using isometry and the correlation structure of $B$, we obtain

$$\mathbb{E}\left[M_tM_{t+h}-(\mathbb{E}[M_t])^2\right]$$

$$=\mathbb{E}\left[\begin{array}{l}\int_{R=0}^{R=+\infty}\int_{s=-\infty}^{s=t}e^{-R(t-s)}\left\{Rb(Q_s-\bar{Q})\rho(\mathrm{d}R)\mathrm{d}s+\sigma\sqrt{R}B(\mathrm{d}R,\mathrm{d}s)\right\}\\ \times\int_{P=0}^{P=+\infty}\int_{u=-\infty}^{u=t+h}e^{-P(t+h-s)}\left\{P\mu(Q_u-\bar{Q})\rho(\mathrm{d}P)\mathrm{d}u+\sigma\sqrt{P}B(\mathrm{d}P,\mathrm{d}u)\right\}\end{array}\right]$$

$$=b^2\mathbb{E}\left[\int_{R=0}^{R=+\infty}\int_{s=-\infty}^{s=t}\int_{P=0}^{P=+\infty}\int_{u=-\infty}^{u=t+h}RPe^{-R(t-s)-P(t+h-u)}(Q_s-\bar{Q})(Q_u-\bar{Q})\rho(\mathrm{d}R)\rho(\mathrm{d}P)\mathrm{d}s\mathrm{d}u\right]$$

$$+\sigma^2\mathbb{E}\left[\int_{R=0}^{R=+\infty}\int_{s=-\infty}^{s=t}\int_{P=0}^{P=+\infty}\int_{u=-\infty}^{u=t+h}\sqrt{RP}e^{-R(t-s)-P(t+h-u)}B(\mathrm{d}R,\mathrm{d}s)B(\mathrm{d}P,\mathrm{d}u)\right] \qquad (A8)$$

$$=b^2\int_{R=0}^{R=+\infty}\int_{s=-\infty}^{s=t}\int_{P=0}^{P=+\infty}\int_{u=-\infty}^{u=t+h}RPe^{-R(t-s)-P(t+h-u)}\bar{V}A_Q(|s-u|)\rho(\mathrm{d}R)\rho(\mathrm{d}P)\mathrm{d}s\mathrm{d}u$$

$$+\sigma^2\mathbb{E}\left[\int_{R=0}^{R=+\infty}\int_{s=-\infty}^{s=t}e^{-2R(t-s)-Rh}m_s(\mathrm{d}R)\mathrm{d}s\right]$$

$$=b^2\bar{V}\int_{R=0}^{R=+\infty}\int_{s=-\infty}^{s=t}\int_{P=0}^{P=+\infty}\int_{u=-\infty}^{u=t+h}RPe^{-R(t-s)-P(t+h-u)}A_Q(|s-u|)\rho(\mathrm{d}R)\rho(\mathrm{d}P)\mathrm{d}s\mathrm{d}u$$

$$+\frac{1}{2}\sigma^2\mathbb{E}[M_t]\int_0^{+\infty}e^{-Rh}\rho(\mathrm{d}R)$$

The last multiple integral in Eq. (A8) can be rewritten using Eq. (1) as follows:

$$\int_{R=0}^{R=+\infty}\int_{s=-\infty}^{s=t}\int_{P=0}^{P=+\infty}\int_{u=-\infty}^{u=t+h}RPe^{-R(t-s)-P(t+h-u)}A_Q(|s-u|)\rho(\mathrm{d}R)\rho(\mathrm{d}P)\mathrm{d}s\mathrm{d}u$$

$$=\int_{R=0}^{R=+\infty}\int_{s=-\infty}^{s=t}\int_{P=0}^{P=+\infty}\int_{u=-\infty}^{u=t+h}\int_{r=0}^{r=+\infty}RPe^{-R(t-s)-P(t+h-u)-r|s-u|}\theta(\mathrm{d}r)\rho(\mathrm{d}R)\rho(\mathrm{d}P)\mathrm{d}s\mathrm{d}u \qquad (A9)$$

$$=\int_{R=0}^{R=+\infty}\int_{P=0}^{P=+\infty}\int_{r=0}^{r=+\infty}RPe^{-Ph}\left(\int_{s=-\infty}^{s=t}\int_{u=-\infty}^{u=t+h}e^{-R(t-s)-P(t-u)-r|s-u|}\mathrm{d}s\mathrm{d}u\right)\theta(\mathrm{d}r)\rho(\mathrm{d}R)\rho(\mathrm{d}P)$$

By elementary calculations, we obtain

$$\int_{s=-\infty}^{s=t}\int_{u=-\infty}^{u=t+h} e^{-R(t-s)-P(t-u)-r|s-u|}\mathrm{d}s\mathrm{d}u = \int_{s=0}^{s=+\infty}\int_{u=0}^{u=+\infty} e^{-Rs-Pu-r|s-u|}\mathrm{d}s\mathrm{d}u + \int_{s=0}^{s=+\infty}\int_{u=-h}^{u=0} e^{-Rs-Pu-r|s-u|}\mathrm{d}s\mathrm{d}u \quad (A10)$$

and

$$\int_{s=0}^{s=+\infty}\int_{u=-h}^{u=0} e^{-Rs-Pu-r|s-u|}\mathrm{d}s\mathrm{d}u = \int_{s=0}^{s=+\infty}\int_{u=-h}^{u=0} e^{-Rs-Pu-r(s-u)}\mathrm{d}s\mathrm{d}u = \frac{1}{R+r}\frac{1}{P-r}\left(e^{(P-r)h}-1\right). \quad (A11)$$

Using Eq. (A7), we obtain

$$\int_{R=0}^{R=+\infty}\int_{s=-\infty}^{s=t}\int_{P=0}^{P=+\infty}\int_{u=-\infty}^{u=t+h} RPe^{-R(t-s)-P(t+h-u)}A_Q(|s-u|)\rho(\mathrm{d}R)\rho(\mathrm{d}P)\mathrm{d}s\mathrm{d}u$$
$$= \int_{R=0}^{R=+\infty}\int_{P=0}^{P=+\infty}\int_{r=0}^{r=+\infty} RPe^{-Ph}\left(\frac{1}{P+R}\left(\frac{1}{P+r}+\frac{1}{R+r}\right) + \frac{1}{R+r}\frac{1}{P-r}\left(e^{(P-r)h}-1\right)\right)\theta(\mathrm{d}r)\rho(\mathrm{d}R)\rho(\mathrm{d}P) \quad (A12)$$

Combining Eqs. (A8)-(A12) along with the following Eq. (A13) yields Eq. (11):

$$\mathbb{E}\left[M_t M_{t+h} - (\mathbb{E}[M_t])^2\right] = \mathbb{V}[M_t]A_M(h). \quad (A13)$$

□

*Proof of Proposition 3*

We derive $F$ and $G$ separately. To obtain $F$, by its definition we have

$$\begin{aligned}
F(l) &= \mathbb{E}\left[(Q_t - \mathbb{E}[Q_t])(M_{t+l} - \mathbb{E}[M_t])\right] \\
&= \mathbb{E}\left[\int_{R=0}^{R=+\infty}\int_{s=-\infty}^{s=t+l} e^{-R(t+l-s)}\left\{Rb(Q_s - \bar{Q})\rho(\mathrm{d}R)\mathrm{d}s + \sigma\sqrt{R}B(\mathrm{d}R,\mathrm{d}s)\right\} \times (Q_t - \bar{Q})\right] \\
&= b\mathbb{E}\left[\int_{R=0}^{R=+\infty}\int_{s=-\infty}^{s=t+l} e^{-R(t+l-s)}R(Q_s - \bar{Q})(Q_t - \bar{Q})\rho(\mathrm{d}R)\mathrm{d}s\right] \\
&= b\int_{R=0}^{R=+\infty}\int_{s=-\infty}^{s=t+l} e^{-R(t+l-s)}R\mathbb{E}\left[(Q_s - \bar{Q})(Q_t - \bar{Q})\right]\rho(\mathrm{d}R)\mathrm{d}s \\
&= b\bar{V}\int_{R=0}^{R=+\infty}\int_{s=-\infty}^{s=t+l} e^{-R(t+l-s)}RA_Q(|t-s|)\rho(\mathrm{d}R)\mathrm{d}s \\
&= b\bar{V}\int_{R=0}^{R=+\infty}\int_{r=0}^{r=+\infty} R\left(\int_{s=-\infty}^{s=t+l} e^{-R(t+l-s)-r|t-s|}\mathrm{d}s\right)\theta(\mathrm{d}r)\rho(\mathrm{d}R) \\
&= b\bar{V}\int_{R=0}^{R=+\infty}\int_{r=0}^{r=+\infty}\left\{\frac{R}{R+r}e^{-Rl} + \frac{R}{R-r}\left(e^{-rl}-e^{-Rl}\right)\right\}\theta(\mathrm{d}r)\rho(\mathrm{d}R)
\end{aligned} \quad (A14)$$

which yields Eq. (16). Here, we used the following elementary calculation for $r \ne R$:

$$\begin{aligned}
\int_{s=-\infty}^{s=t+l} e^{-R(t+l-s)-r|t-s|}\mathrm{d}s &= \int_{-\infty}^{t} e^{-R(t+l-s)-r|t-s|}\mathrm{d}s + \int_{t}^{t+l} e^{-R(t+l-s)-r|t-s|}\mathrm{d}s \\
&= \int_{-\infty}^{t} e^{-R(t+l-s)-r(t-s)}\mathrm{d}s + \int_{t}^{t+l} e^{-R(t+l-s)-r(s-t)}\mathrm{d}s \\
&= e^{-Rl}\int_{-\infty}^{t} e^{-(R+r)(t-s)}\mathrm{d}s + e^{-Rl}\int_{t}^{t+l} e^{-(R-r)(t-s)}\mathrm{d}s \\
&= \frac{1}{R+r}e^{-Rl} + \frac{1}{R-r}\left(e^{-rl}-e^{-Rl}\right)
\end{aligned} \quad (A15)$$

The case $r = R$ is understood through the classical l'Hôpital's rule.

To obtain $G$, by its definition we have

$$\begin{aligned}
G(l) &= \mathbb{E}\left[(Q_t - \mathbb{E}[Q_t])(M_{t-l} - \mathbb{E}[M_t])\right] \\
&= \mathbb{E}\left[(Q_{t+l} - \bar{Q})(M_t - \mathbb{E}[M_t])\right] \\
&= \mathbb{E}\left[\int_{R=0}^{R=+\infty}\int_{s=-\infty}^{s=t} e^{-R(t-s)}\left\{Rb(Q_s - \bar{Q})\rho(\mathrm{d}R)\mathrm{d}s + \sigma\sqrt{R}B(\mathrm{d}R,\mathrm{d}s)\right\} \times (Q_{t+l} - \bar{Q})\right] \\
&= b\mathbb{E}\left[\int_{R=0}^{R=+\infty}\int_{s=-\infty}^{s=t} e^{-R(t-s)}R(Q_s - \bar{Q})(Q_{t+l} - \bar{Q})\rho(\mathrm{d}R)\mathrm{d}s\right] \\
&= b\int_{R=0}^{R=+\infty}\int_{s=-\infty}^{s=t+l} e^{-R(t-s)}R\mathbb{E}\left[(Q_s - \bar{Q})(Q_{t+l} - \bar{Q})\right]\rho(\mathrm{d}R)\mathrm{d}s \\
&= b\bar{V}\int_{R=0}^{R=+\infty}\int_{s=-\infty}^{s=t+l} e^{-R(t-s)}RA_Q(t+l-s)\rho(\mathrm{d}R)\mathrm{d}s \\
&= b\bar{V}\int_{R=0}^{R=+\infty}\int_{r=0}^{r=+\infty} R\left(\int_{s=-\infty}^{s=t+l} e^{-R(t-s)-r(t+l-s)}\mathrm{d}s\right)\theta(\mathrm{d}r)\rho(\mathrm{d}R) \\
&= b\bar{V}\int_{R=0}^{R=+\infty}\int_{r=0}^{r=+\infty} \frac{R}{R+r}e^{-rl}\theta(\mathrm{d}r)\rho(\mathrm{d}R)
\end{aligned} \quad (A16)$$

which yields Eq. (17). Here, we used $t + l \geq s$. The representation (18) follows by setting $h = 0$ in $F(h)$ or $G(h)$.

Finally, the right-hand sides of Eqs. (16) and (17) exist if $\int_0^{+\infty} R\rho(\mathrm{d}R) < +\infty$ because of

$$\begin{aligned}
F(l) &= b\bar{V}\int_{R=0}^{R=+\infty}\int_{r=0}^{r=+\infty}\left\{\frac{R}{R+r}e^{-Rl} + \frac{R}{R-r}(e^{-rl} - e^{-Rl})\right\}\theta(\mathrm{d}r)\rho(\mathrm{d}R) \\
&\leq b\bar{V}\int_{R=0}^{R=+\infty}\int_{r=0}^{r=+\infty}\left\{e^{-Rl} + lR\max\{e^{-Rl}, e^{-rl}\}\right\}\theta(\mathrm{d}r)\rho(\mathrm{d}R) \\
&\leq b\bar{V}\int_{R=0}^{R=+\infty}\int_{r=0}^{r=+\infty}\left\{e^{-Rl} + lR\right\}\theta(\mathrm{d}r)\rho(\mathrm{d}R) \qquad , l \geq 0 \\
&\leq b\bar{V}\left(1 + l\int_0^{+\infty} R\rho(\mathrm{d}R)\right) \\
&< +\infty
\end{aligned} \quad (A17)$$

and

$$G(l) \leq b\bar{V}\int_{R=0}^{R=+\infty}\int_{r=0}^{r=+\infty} e^{-rl}\theta(\mathrm{d}r)\rho(\mathrm{d}R) = b\bar{V}A_Q(l) < +\infty, \; l \geq 0. \quad (A18)$$

Here, we used the elementary calculation $\frac{1}{R-r}(e^{-rl} - e^{-Rl}) \leq l\max\{e^{-Rl}, e^{-rl}\}$ for $r, R > 0$, where the case $r = R$ is understood through the classical l'Hôpital's rule.

□

## A2. Finite-dimensional version

A finite-dimensional counterpart of the proposed model is presented, along with its associated generalized Riccati equation. This setting is analogous to that discussed in **Section 3.3**. A repetitive discussion is provided here for clarity. In this section, we assume the superposition of the Ornstein–Uhlenbeck processes for discharge $Q$ (**Section 4.2**).

First, we discretize the probability measures $\rho$ and $\pi$ under the assumption that $\rho$ admits a probability density:

$$\rho(\mathrm{d}R) \approx \frac{1}{N}\sum_{i=1}^{N}\delta(R-R_i) \text{ and } \pi(\mathrm{d}r) \approx \frac{1}{N}\sum_{j=1}^{N}\delta(R-r_j). \tag{A19}$$

Here, $N \in \mathbb{N}$ is the degree of freedom of discretization and $R_i$ is the $\frac{2i-1}{2N}$ th quantile level of the probability measure $\rho$ : $\int_0^{R_i} \rho(\mathrm{d}R) = \frac{2i-1}{2N}$. Similarly, $r_j$ is the $\frac{2j-1}{2N}$ th quantile level: $\int_0^{r_j} \pi(\mathrm{d}r) = \frac{2j-1}{2N}$. Accordingly, we set the finite-dimensional versions of the memory process and discharge:

$$M_t^{(N)} = \sum_{i=1}^{N} m_t^{(N,i)} \text{ and } Q_t^{(N)} = \sum_{j=1}^{N} q_t^{(N,j)}, \; t \in \mathbb{R} \tag{A20}$$

with

$$\mathrm{d}m_t^{(N,i)} = -R_i\left(m_t^{(N,i)} - \frac{1}{N}\left(a + bQ_t^{(N)}\right)\right)\mathrm{d}t + \sigma\sqrt{R_i m_t^{(N,i)}}\,\mathrm{d}W_t^i, \; i = 1,2,3,...,N \tag{A21}$$

and

$$\mathrm{d}q_t^{(N,j)} = -r_j q_t^{(N,j)}\mathrm{d}t + z^{(N,j)}\mathrm{d}\Xi_t^j, \; j = 1,2,3,...,N. \tag{A22}$$

Here, each $W_t^i$ is an independent one-dimensional standard Brownian motion that is mutually independent of others. The discretized discharge $Q$ and each $\Xi_t^j$ is an independent Poisson random measure of intensity for each $N^{-1}\mathrm{d}s\nu\left(\mathrm{d}z^{(N,j)}\right)$ of the tempered stable model, with $z^{(N,j)}$ denoting jump size. Each $W_t^i$ and $\Xi_t^j$ are mutually independent.

The SDEs (A21) and (A22), as a finite-dimensional version of the proposed model, are affine (e.g., Duffie et al., 2003)[57]. Thus, the moment-generating function $\mathbb{E}\left[e^{\varpi M_t^{(N)}}\right]$ ($\varpi \leq 0$), at a stationary state, can be obtained by solving a generalized Riccati equation. For a generic smooth function $\phi = \phi(m_1, m_2, ..., m_N, q_1, q_2, ..., q_N)$ with $(m_1, m_2, ..., m_N, q_1, q_2, ..., q_N) \in \mathbb{R}^{2N}$, the infinitesimal generator associated with Eqs. (A21) and (A22) is obtained as

$$G\phi = \sum_{i=1}^{N}\left\{-R_i\left(m_i - \frac{1}{N}\left(a + b\sum_{j=1}^{N}q_j\right)\right)\frac{\partial\phi}{\partial m_i} + \frac{1}{2}\sigma^2 R_i m_i \frac{\partial^2\phi}{\partial m_i^2}\right\}$$
$$+ \sum_{j=1}^{N}\left\{-r_j q_j \frac{\partial\phi}{\partial q_j} + \frac{1}{N}\int_0^{+\infty}\Delta\phi(z_j)\nu(\mathrm{d}z_j)\right\} \tag{A23}$$

with

$$\Delta\phi(z_j) = \phi(m_1, m_2, ..., m_N, q_1, q_2, ..., q_j + z_j, ..., q_N) - \phi. \tag{A24}$$

The right-hand side term of Eq. (A23) denotes a partial differential operation with differential operators multiplied by the affine coefficients for each $m_i, q_j$, which is why the proposed model is categorized as an affine process. The affine structure is important in applications. For example, the affine structure is the key to obtaining the generalized Riccati equation associated with Eqs. (A21) and (A22). We assume the exponential functional form to estimate the form of the conditional moment-generating function (conditioned at time 0):

$$\phi(t, m_1, m_2, ..., m_N, q_1, q_2, ..., q_N) = \mathbb{E}\left[ e^{\varpi M_t^{(N)}} \Big| (m_1, m_2, ..., m_N, q_1, q_2, ..., q_N) \right]$$
$$= \exp\left( \kappa_t + \sum_{i=1}^{N} \omega_t^{(i)} m_i + \sum_{j=1}^{N} \psi_t^{(j)} q_j \right), \quad (A25)$$

which is substituted in the Kolmogorov equation:

$$\frac{\partial \Phi}{\partial t} = G\Phi, \quad t > 0, \quad (m_1, m_2, ..., m_N, q_1, q_2, ..., q_N) \in \mathbb{R}^{2N}. \quad (A26)$$

Then, we obtain an identity to determine each coefficient $\kappa_t, \omega_t^{(\cdot)}, \psi_t^{(\cdot)}$:

$$\frac{d}{dt}\left( \kappa_t + \sum_{i=1}^{N} \omega_t^{(i)} m_i + \sum_{j=1}^{N} \psi_t^{(j)} q_j \right)$$
$$= \sum_{i=1}^{N} \left\{ -R_i \left( m_i - \frac{1}{N}\left( a + b\sum_{j=1}^{N} q_j \right) \right) \omega_t^{(i)} + \frac{1}{2}\sigma^2 R_i m_i \left( \omega_t^{(i)} \right)^2 \right\}$$
$$+ \sum_{j=1}^{N} \left\{ -r_j q_j \psi_t^{(j)} + \frac{1}{N}\int_0^{+\infty} \left( \exp\left( \psi_t^{(j)} z_j \right) - 1 \right) \nu(dz_j) \right\} \quad . \quad (A27)$$
$$= \sum_{i=1}^{N} m_i \left\{ -R_i \omega_t^{(i)} + \frac{1}{2}\sigma^2 R_i \left( \omega_t^{(i)} \right)^2 \right\}$$
$$+ \sum_{j=1}^{N} q_j \left( -r_j \psi_t^{(j)} + b\frac{1}{N}\sum_{i=1}^{N} R_i \omega_t^{(i)} \right) + \frac{1}{N}\sum_{i=1}^{N} a R_i \omega_t^{(i)} + \frac{1}{N}\sum_{j=1}^{N} \int_0^{+\infty} \left( \exp\left( \psi_t^{(j)} z_j \right) - 1 \right) \nu(dz_j)$$

Consequently, we derive the generalized Riccati equation as follows:

$$\frac{d\kappa_t}{dt} = \frac{1}{N} a \sum_{i=1}^{N} R_i \omega_t^{(i)} + \frac{1}{N} \sum_{j=1}^{N} \int_0^{+\infty} \left( \exp\left( \psi_t^{(j)} z_j \right) - 1 \right) \nu(dz_j), \quad t > 0, \quad (A28)$$

$$\frac{d\omega_t^{(i)}}{dt} = -R_i \omega_t^{(i)} + \frac{1}{2}\sigma^2 R_i \left( \omega_t^{(i)} \right)^2, \quad 1 \le i \le N, \quad t > 0, \quad (A29)$$

$$\frac{d\psi_t^{(j)}}{dt} = -r_j \psi_t^{(j)} + b\frac{1}{N}\sum_{i=1}^{N} R_i \omega_t^{(i)}, \quad 1 \le j \le N, \quad t > 0 \quad (A30)$$

subject to the initial conditions $\omega_0^{(\cdot)} = \varpi$ and $\psi_0^{(\cdot)} = \kappa_0 = 0$. The generalized Riccati equation associated with the proposed model is obtained by formally taking $N \to +\infty$ using the notation:

$$\frac{d\kappa_t}{dt} = a\int_0^{+\infty} R\omega_t(R) \rho(dR) + a\int_0^{+\infty} \int_0^{+\infty} \left( \exp\left( \psi_t(r) z \right) - 1 \right) \nu(dz) \pi(dr), \quad t > 0, \quad (A31)$$

$$\frac{d\omega_t(R)}{dt} = -R\omega_t(R) + \frac{1}{2}\sigma^2 R\left( \omega_t(R) \right)^2, \quad R > 0, \quad t > 0, \quad (A32)$$

$$\frac{d\psi_t(r)}{dt} = -r\psi_t(r) + b\int_0^{+\infty} R\omega_t(R) \rho(dR), \quad r > 0, \quad t > 0 \quad (A33)$$

subject to the initial conditions of $\omega_0(0) = \varpi$ and $\psi_0(0) = \kappa_0 = 0$. These equations, which is a system of non-local ordinary differential equations, are mathematically interesting by itself but are not directly studied in detail in this study because it is not used in our analysis. This point is currently undergoing in a separate study from mathematical viewpoints.

**A3. Computation of the proposed model**

The proposed model is numerically discretized in **Section 4** to compute the probability densities and statistics considering three water quality indices. This section briefly discusses the employed computational resolution. **Tables A1** to **A3** present the computed average and variance for TN, TP, and TOC against the three different computational resolutions $(\Delta t, N)$: coarse $(\Delta t, N) = (0.02, 1024)$, nominal $(\Delta t, N) = (0.01, 2048)$, and fine $(\Delta t, N) = (0.005, 4096)$. The sampling period is 1,000 years, with a 1,000-year burn-in period. **Fig. A1** presents a comparison of empirical and computed PDFs for TN (mg/L) considering different computational resolutions. The computational results overlap in the ordinary scale and have similar tails in the logarithmic scale with less scattering near the tail as the resolution increases. As presented in **Tables A1** to **A3**, increasing the resolution does not necessarily improve the computational accuracy of statistics as those for fine resolution because increasing the degree of freedom $N$, numerically deals with SDEs considering smaller $r, R$, requiring a longer burn-in period for computing statistics that are sufficiently close to be stationary, particularly in computing high-order statistics. Accordingly, a nominal resolution presents a reasonably balanced choice. Another probable reason is that the parameter estimation is conducted considering the resolution $N = 2,048$. Thus, we analyze model uncertainty, as discussed in **Section 4**.

**Table A1** Computed statistics and different computational resolutions for TN

|  | Theoretical | Coarse | Nominal | Fine |
|---|---|---|---|---|
| Average | 1.084.E+00 | 1.074.E+00 | 1.074.E+00 | 1.085.E+00 |
| Relative error |  | 9.428.E-03 | 9.031.E-03 | 4.689.E-04 |
| Variance | 2.676.E-01 | 2.601.E-01 | 2.703.E-01 | 2.548.E-01 |
| Relative error |  | 2.808.E-02 | 1.015.E-02 | 4.792.E-02 |

**Table A2** Computed statistics and different computational resolutions for TP

|  | Theoretical | Coarse | Nominal | Fine |
|---|---|---|---|---|
| Average | 1.248.E+00 | 1.218.E+00 | 1.223.E+00 | 1.244.E+00 |
| Relative error |  | 2.405.E-02 | 1.986.E-02 | 2.910.E-03 |
| Variance | 1.659.E+00 | 1.591.E+00 | 1.651.E+00 | 1.853.E+00 |
| Relative error |  | 4.087.E-02 | 4.909.E-03 | 1.171.E-01 |

**Table A3** Computed statistics and different computational resolutions for TOC

|  | Theoretical | Coarse | Nominal | Fine |
|---|---|---|---|---|
| Average | 1.119.E+00 | 1.084.E+00 | 1.104.E+00 | 1.114.E+00 |
| Relative error |  | 3.143.E-02 | 1.364.E-02 | 4.244.E-03 |
| Variance | 5.432.E-01 | 4.919.E-01 | 5.562.E-01 | 6.497.E-01 |
| Relative error |  | 9.453.E-02 | 2.377.E-02 | 1.960.E-01 |

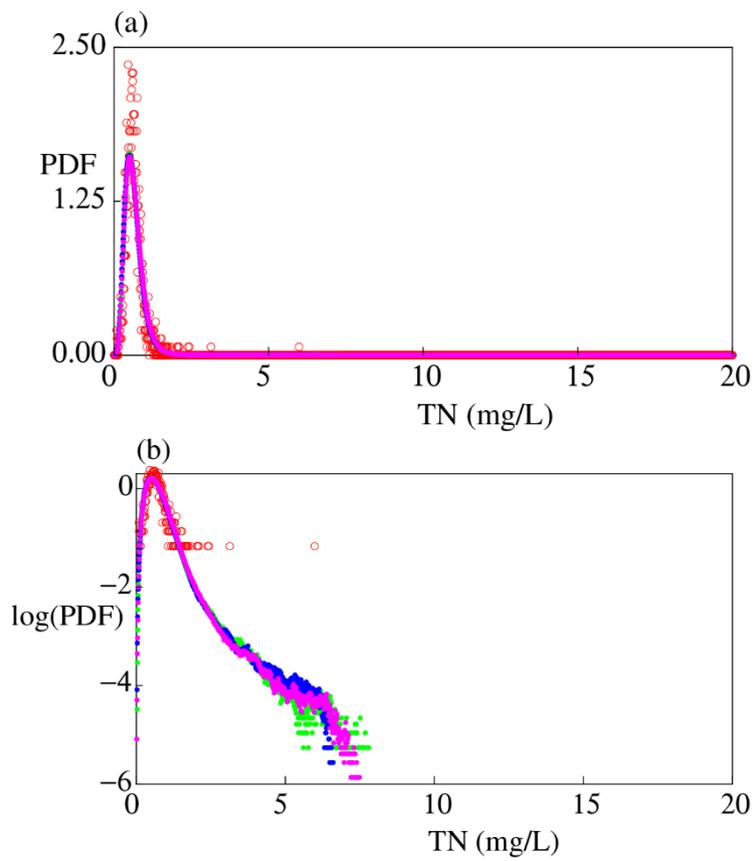

**Fig. A1** Comparison of empirical and computed PDFs for TN on (a) ordinary scale and (b) ordinary logarithmic scale: empirical results (red), and computational ones with coarse (green), nominal (blue), and resolutions (magenta)